\documentclass[11pt]{article}
\usepackage{amsmath,amsfonts,amssymb}
\usepackage{multirow} 
\usepackage{float}

\usepackage[section]{placeins} 

\usepackage{graphicx}
\usepackage{epstopdf}
\usepackage{bm}
\usepackage{color}
\usepackage{lscape}
\usepackage{url}
\usepackage{hyperref}

\DeclareMathOperator\erfc{erfc}
\DeclareMathOperator\erf{erf}

\newcommand{\bd}[1]{\mathbf{#1}}

\newcommand{\tb}{\mathbf}
\newcommand{\ff}{{\tb f}}

\newcommand{\uu}{{\tb u}}
\newcommand{\vv}{{\tb v}}
\newcommand{\ww}{{\tb w}}
\newcommand{\xx}{{\tb x}}
\newcommand{\yy}{{\tb y}}
\newcommand{\zz}{{\tb z}}
\newcommand{\rr}{{\tb r}}
\newcommand{\nn}{{\tb n}}

\newcommand{\pp}{{\tb p}}
\newcommand{\UU}{{\tb U}}

\newcommand{\gbar}{\bar{g}}

\newcommand{\xxh}{{\hat \xx}}

\newcommand{\del}{\delta}

\newcommand{\lam}{\lambda}
\newcommand{\sig}{\sigma}
\newcommand{\kap}{\kappa}
\newcommand{\kapz}{\kappa_0}
\newcommand{\pa}{\partial}
\newcommand{\beq}{\begin{equation}}
\newcommand{\eeq}{\end{equation}}
\newcommand{\Sl}{\mathcal{S}}
\newcommand{\Dl}{\mathcal{D}}

\definecolor{red}{rgb}{1.0,0.0,0.0}

\title{
High Order Regularization of \\ Nearly Singular Surface Integrals}

\author{J. Thomas Beale\thanks{Department of Mathematics, Duke University, 
Durham, NC, 27708 USA beale@math.duke.edu} 
\and Svetlana Tlupova\thanks{Department of Mathematics,
Farmingdale State College, SUNY, Farmingdale, NY 11735, USA tlupovs@farmingdale.edu}}

\date{\today}

\begin{document}

\maketitle

\begin{abstract}

Solutions of partial differential equations can often be written as surface integrals having a kernel related to a singular fundamental solution.  Special methods are needed to evaluate the integral accurately at points on or near the surface.  Here we derive formulas to regularize the integrals with high accuracy, using analysis from~\cite{extrap}, so that a standard quadrature can be used without special care near the singularity.
We treat single or double layer integrals for harmonic functions or for Stokes flow.  The nearly singular case, evaluation at points close to the surface, can be needed when surfaces are close to each other, or to find values at grid points near a surface.  We derive formulas for regularized kernels with error $O(\delta^p)$ where $\delta$ is the smoothing radius and $p = 3$, $5$, $7$.  With spacing $h$ in the quadrature, we choose $\delta = \kappa h^q$ with $q<1$ so that the discretization error is controlled as $h \to 0$.  We see the predicted order of convergence $O(h^{pq})$ in various examples. 
Values at all grid points can be obtained from those near the surface in an efficient manner suggested in~\cite{mayo85}.  With this technique we obtain high order accurate grid values for a harmonic function determined by interfacial conditions and for the pressure and velocity in Stokes flow around a translating spheroid.

\end{abstract}



\section{Introduction}

Solutions of partial differential equations, determined by boundary or interfacial conditions, can often be written as surface integrals having a kernel related to a singular fundamental solution.  Special methods are needed to evaluate the integral accurately at points on or near the surface.  Here we derive formulas to regularize the integrals with high accuracy, using analysis from~\cite{extrap}, so that a standard quadrature can be used without special care near the singularity.
We treat single or double layer integrals for harmonic functions or for Stokes flow.  The singular case, evaluation on the surface, is needed for solving integral equations.  The nearly singular case, evaluation at points close to the surface, can be needed when surfaces are close to each other, or to find values at grid points near a surface.  Values at all grid points can be obtained from those near the surface in an efficient manner suggested in~\cite{mayo85}.  Formulas derived here are valid for evaluation on the surface as well as at nearby points, with simpler expressions given for the special case of points on the surface.

To regularize the integral we replace
a singular kernel such as $1/r$ with a smooth version, in our case
$\erf(r/\delta)/r$.  Here $\delta$ is a numerical parameter which determines the radius of smoothing.  We begin with a simple choice of the smoothing factor, or shape factor, and modify it to achieve a high order error in the integral, $O(\delta^p)$ with $p = 3,5,7$.  The modification is based on analytic expressions for the error, when the integral is evaluated near a point on the surface, derived in~\cite{extrap}; see~\eqref{basic} below.  This higher order permits $\delta$ to be larger, in order to control the error in discretizing the regularized integral.  It is natural to choose $\delta$ proportional to $h$, the mesh spacing in the quadrature, and with this choice we see total errors of $O(h^p)$ in examples.  However, if 
$\delta/h$ is fixed as $h$ is reduced, we expect low order accuracy in $h$ in the quadrature.  Thus to obtain convergence as $h \to 0$,
we need to increase $\delta/h$.  We do this by setting $\delta = \kap h^q$ with $q < 1$.
For example with $p = 7$ and $q = 5/7$, we see convergence to $O(h^5)$.  In~\cite{extrap} we used extrapolation with several choices of $\delta$ to obtain higher order, rather than modifying the smoothing factor.  The present work extends and applies the results of~\cite{extrap}.

To compute the regularized integrals we need to choose a quadrature rule. The total error consists of the smoothing error plus the error in discretizing the regularized integral.
In this work we use a quadrature rule due to J. Wilson~\cite{wilson-10},~\cite{byw}.  It uses a partition of unity on the unit sphere, which is applied to the normal vector at points on the surface; this is different from the usual role of partitions of unity.  The quadrature is reduced to sums in the coordinate planes, cut off by the partition functions.  For smooth integrands and surfaces it has the high order accuracy of the trapezoidal rule without boundaries.  For our regularized integrals the accuracy depends on the balance between $\delta$ and the mesh size $h$.  This rule and its accuracy are discussed further in Sect. 4 of~\cite{extrap}.  It has the advantages that only limited knowledge of the surface is needed, and the points are regularly spaced in coordinate planes.

A variety of numerical examples are presented in Sect. 4 using high order smoothing.
We verify the predicted order of accuracy in simple examples at grid points within distance $h$ of the surface,
with $\delta$ proportional to $h$ or $h^q$.  For
harmonic functions with specified jump conditions at an interface, we
use the integral representation to obtain the solution at arbitrary grid points, using a version of the
procedure of~\cite{mayo85}:  We compute the solution at grid points near the interface as nearly singular
integrals with $O(h^5)$ accuracy.  We form a fourth order discrete Laplacian, extend it by zero away from the interface, and invert using the discrete Fourier sine transform, to obtain the solution on the entire grid.  The solution and its computed gradient are accurate to $O(h^4)$.  The reasons for this accuracy are discussed in Sect. 4.
In a similar manner we compute the Stokes flow around a translating spheroid, first computing the pressure $p$ and $\nabla p$ and then the velocity and its gradient.  In another example we solve an integral equation for the velocity in Stokes flow with two spheroids close to each other.

A large amount of work has been done on computational methods for singular integrals.
A portion of this work has concerned nearly singular integrals on surfaces.  Often
values close to the surface are obtained by extrapolating from interior values, e.g. in~\cite{ying-biros-zorin-06}, and in quadrature by expansion (QBX) or hedgehog methods
 \cite{klinttorn,klockner-barnett-greengard-oneil-13,zorincplx,siegel-tornberg-18}. 
With the singularity subtraction technique \cite{helsing-13,planewave} a most singular part is evaluated
analytically leaving a more regular remainder. In \cite{perez}, and recently in \cite{jiangzhu},
 a harmonic approximation to the density function is used to reduce the singularity.
In \cite{nitsche} and more recently \cite{nitschewu}, local corrections are used
that are determined from analytic expressions.
Regularization of the sort used here has been applied extensively to problems in biology
modeled by Stokes flow \cite{cortez2d,cortez3d,cortezsurf}.
Methods based on heat potentials,
e.g. \cite{lightweight}, use asymptotic analysis related to that in \cite{extrap}.

To describe our approach, we outline the treatment of the single layer integral
\beq \label{sgllayer}
  \Sl(\yy) = \int_\Gamma G(\xx - \yy)f(\xx)\,dS(\xx)\,, \quad
  G(\rr)  =  -\frac{1}{4\pi|\rr|} \eeq
on a closed surface $\Gamma$ with given density function $f$.  To begin we define the regularized version
\beq \label{sglreg}
   \Sl_\del(\yy) = \int_\Gamma G_\del(\xx - \yy)f(\xx)\,dS(\xx)\,, \quad     G_\delta(\rr) = G(\rr)s_1(|\rr|/\del) \eeq
with 
\beq \label{sglshape}
s_1(\rho) = \erf(\rho) = \frac{2}{\sqrt{\pi}}\int_0^\rho e^{-\sig^2}\,d\sig \eeq
Then $G_\delta$ is a smooth function of $\xx-\yy$ and
$\erf(\rho) \to 1$ rapidly as $\rho$ increases.  
The singularity at $\yy = \xx$ is replaced by a peak of $O(1/\del)$.
If $\yy$ is near $\Gamma$, then $\Sl_\del - \Sl = O(\del)$.
For such $\yy$ we can write 
\beq \label{closest} \yy = \xx_0 + b\nn \eeq
 where $\xx_0$ is the closest point to $\yy$
on $\Gamma$, $\nn$ is the outward normal vector at $\xx_0$, and $b$ is the signed distance.
We wish to improve the error by replacing $s_1(\rho)$ with a version giving higher order.
In \cite{extrap}, Sect. 3, we showed that
\beq \label{basic}
 \Sl_\del(\yy) = \Sl(\yy) +  C_1\del I_0(b/\del)  + C_2\del^3 I_2(b/\del)
         + C_3\del^5 I_4(b/\del)  + O(\del^7)  \eeq 
uniformly for $\yy$ near the surface.  Here $I_0$, $I_2$, $I_4$ are certain integrals that appear in the analysis and are known explicitly; see \eqref{eyes}--\eqref{eye4} below.  The coefficients $C_1$, $C_2$, $C_3$ depend on the specifics of the problem as well as $\yy$, but not on $\del$, and are difficult to find.  Our strategy here is to modify $s_1$ in such a way that the integrals $I_n$ in \eqref{basic} are replaced by zero, so that the leading error is eliminated.   We do this in a systematic manner so that a similar procedure can be used for the other integral kernels, and so that the expressions obtained can easily be truncated to lower order as desired.  The modified version of $s_1$ derived in Sect. 2 has the form
(see \eqref{s17})
\beq \label{s17g} s_1^{(7)}(\rho) = \erf(\rho) + \left(c_1\rho + c_2\rho^3 + c_3\rho^5\right)e^{-\rho^2} \eeq
with $c_1$, $c_2$, $c_3$ depending on $b/\del$. If we replace $s_1$ in \eqref{sglreg} with $s_1^{(7)}$,
the corresponding integral $\Sl_\del^{(7)}$ has error reduced to $O(\del^7)$.  
Lower order versions of $s_1$ have accuracy $O(\del^3)$ or $O(\del^5)$ with fewer terms. 
Specific formulas are given in Sect. 2, as well as similar formulas for the other integral kernels.
Since $1 - \erf{\rho}$ and $\exp(-\rho^2)$ decay rapidly as $\rho$ increases, the smoothing factor
can be neglected for $\rho = r/\del$ large enough, e.g. $\rho > 8$ or $r > 8\del$.

For efficiency of evaluating the discretized integrals, fast summation methods suitable for regularized kernels~\cite{rbfstokes, wang-krasny-tlupova-20, yingradial} could be used. Recently in~\cite{siebor-tlupova}, the kernel-independent treecode of~\cite{wang-krasny-tlupova-20} is used to compute the regularized integrals with about $O(N\log N)$ efficiency, where $N$ is the system size. The regularization is localized and the far-field particle-particle interactions are evaluated efficiently through particle-cluster interactions. Experiments carried out in~\cite{siebor-tlupova} determine optimal treecode parameters in the mesh size $h$, both in terms of accuracy and efficiency of the computations.

Formulas are presented in Sect. 2 for regularization of various nearly singular integrals, including
single and double layer potentials for harmonic functions, single layer, or Stokeslet, integrals for
Stokes flow, and double layer, or stresslet, integrals for Stokes flow.  These are used for evaluation
near the surface.  Sect. 3 gives simplified formulas for the special case of evaluation on the surface.  Numerical examples are presented in Sect. 4. An appendix has a formula for extending a smooth function on the boundary of a computational box to the interior. The source code for numerical examples is available for download (github.com/stlupova/high-order-regularization).


\section{Formulas for points near the surface}
\label{sec:offSurf}

We first derive the versions of the smoothing factor in the form \eqref{s17g} 
to be used with \eqref{sglreg} in evaluating the single layer potential
and then proceed to the other cases.  We assume throughout that the point $\yy$ is close to the surface,
so that \eqref{closest} applies.
Let $\lam = b/\delta$.
The integrals occurring in \eqref{basic} are
\beq \label{eyes} I_n(\lam) = \int_0^\infty 
    \frac{\erfc(\sqrt{\sig^2 + \lam^2})}
       {\sqrt{\sig^2 + \lam^2}} \sig^{n+1}\,d\sig 
   = \int_{|\lam|}^\infty \erfc(\rho) (\rho^2 - \lam^2)^{n/2}\,d\rho \,, \quad n = 0,2,4
         \eeq
where $\erfc$ is the complementary error function, $\erfc = 1 - \erf$.  They are
\beq \label{eye0}
I_0(\lambda) = e^{-\lambda^2}/\sqrt{\pi}
                         - |\lambda|\erfc|\lambda|  \eeq
\beq \label{eye2}
I_2(\lambda) = \frac23\left( 
    (\frac12 - \lambda^2) e^{-\lambda^2}/\sqrt{\pi}
       + |\lambda|^3 \erfc|\lambda| \right) \eeq
\beq \label{eye4}
I_4(\lambda) =  \frac{8}{15}\left(
   (\frac34 - \frac12\lam^2 + \lam^4)e^{-\lambda^2}/\sqrt{\pi}
- |\lambda|^5 \erfc|\lambda| \right) \eeq
The integrals $I_n$ could be thought of as generalized moments for $(1 - s_1)/r = \erfc/r$.

\medskip

{ \bf The single layer potential.}
We want to modify $s_1$ to $s_1^{(7)}$ so that the corresponding quantities in \eqref{basic} are replaced by zero.
We seek $s_1^{(7)}$ in the form
\beq 
	s_1^{(7)}(\rho) = s_1(\rho) + a_1 \rho s_1'(\rho) + a_2 \rho^2 s_1''(\rho) + a_3 \rho^3 s_1'''(\rho) 
\eeq
with $' = d/d\rho$ and coefficients $a_1$, $a_2$, $a_3$ to be determined. The integral in place of~\eqref{eyes} will have extra terms of the form
\beq I_{kn}(\lambda) = \int_{|\lam|}^\infty \rho^k s_1^{[k]}(\rho) (\rho^2 - \lam^2)^{n/2}\,d\rho\,, 
       \quad k= 1,2,3; \; n = 0,2,4 \eeq
and our requirement becomes
\beq \label{sys1} a_1 I_{1n} + a_2 I_{2n} + a_3 I_{3n} = I_n \,,\quad n=0,2,4 \eeq
a system of linear equations.  The integrals $I_{kn}$ are easily found, 
and the solution of the system is
\beq \label{tri1}
 a_3 = \frac{\sqrt{\pi}}{16} (2I_0 - 4I_2 + I_4) \, e^{\lam^2}  \eeq
 \beq \label{tri2}
 a_2 = \frac{\sqrt{\pi}}{2} (I_0 - I_2) \, e^{\lam^2}  + (4\lam^2+7) \, a_3 \eeq
 \beq \label{tri3}
 a_1 = \sqrt{\pi} I_0 \, e^{\lam^2}  + 2(\lam^2+1)\, a_2 - (4\lam^4+6\lam^2+6)\, a_3 \eeq
Using expressions for the derivatives $s_1^{[k]}$ we obtain the formula
\beq
	\label{s17} 
	s_1^{(7)}(\rho) = \erf(\rho) + \frac{2}{\sqrt{\pi}}  \left(
   a_1\rho - 2(a_2+a_3)\rho^3 + 4a_3\rho^5 \right)e^{-\rho^2}  
\eeq 
thereby identifying $c_1$, $c_2$, $c_3$ in \eqref{s17g}.  We have written \eqref{tri1}--\eqref{tri3} in
triangular form so that lower order expressions are easily obtained:  The fifth order
version $s_1^{(5)}$ is found by setting $a_3 = 0$ and solving \eqref{tri2}, \eqref{tri3}
for $a_2^{(5)}$ and $a_1^{(5)}$, to obtain
\beq \label{a5} 
a_2^{(5)} = \frac{\sqrt{\pi}}{2} (I_0 - I_2) \, e^{\lam^2}\,, \quad
a_1^{(5)} = \sqrt{\pi} I_0 \, e^{\lam^2}  + 2(\lam^2+1)\, a_2^{(5)}  \eeq
In place of \eqref{s17} we get
\beq  
	\label{s15} 
	s_1^{(5)}(\rho) = \erf(\rho) + \frac{2}{\sqrt{\pi}}  \left(
   a_1^{(5)}\rho - 2a_2^{(5)}\rho^3 \right)e^{-\rho^2} 
\eeq
Similarly, the third order version $s_1^{(3)}$ is obtained
by setting $a_2 = a_3 = 0$ and solving \eqref{tri3} to find $a_1^{(3)} = \sqrt{\pi} I_0 \, e^{\lam^2}$ and
\beq \label{s13}  s_1^{(3)}(\rho) = \erf(\rho) + 2e^{\lam^2} I_0(\lam)  \rho e^{-\rho^2} \eeq

\medskip

{\bf The double layer potential.}
The double layer integral has the form
\beq \label{dbllayer}
  \Dl(\yy) = \int_\Gamma \frac{\pa G(\xx-\yy)}{\pa\nn(\xx)}
        g(\xx)\,dS(\xx) \eeq
We use the familiar identity
\beq \int_\Gamma \frac{\pa G(\xx-\yy)}{\pa\nn(\xx)}\,dS(\xx) = \chi(\yy) \eeq
with $\chi = 1$ for $\yy$ inside, $\chi = 0$ outside, and
$\chi = \frac12$ on $\Gamma$.  We rewrite \eqref{dbllayer} as
\beq \label{dblsub} 
  \Dl(\yy) = \int_\Gamma \frac{\pa G(\xx-\yy)}{\pa\nn(\xx)}
 [g(\xx) - g(\xx_0)]\,dS(\xx) + \chi(\yy)g(\xx_0) \eeq
where again $\xx_0$ is the closest point on $\Gamma$.
The subtraction reduces the singularity and it is necessary for our method.
We now replace $\nabla G$
with the gradient of the smooth function $G_\del$, obtaining
\beq \label{gradreg}
  \nabla G_\del(\rr) = \nabla G(\rr)s_2(|\rr|/\del) = 
    \frac{\rr}{4\pi|\rr|^3}s_2(|\rr|/\del)  \eeq
with
\beq \label{dblshape}
     s_2(\rho) = \erf(\rho) - (2/\sqrt{\pi})\rho e^{-\rho^2}  \eeq
and our regularized form of \eqref{dbllayer} is 
\beq \label{dblreg}
  \Dl_\del(\yy) = \int_\Gamma 
  \frac{\rr\cdot\nn(\xx)}{4\pi|\rr|^3}s_2(|\rr|/\del)
  [g(\xx) - g(\xx_0)]\,dS(\xx) + \chi(\yy)g(\xx_0)\,,
             \quad    \rr = \xx - \yy  \eeq
             
In \cite{extrap} we showed that 
\beq \label{dblerr} \Dl_\del = \Dl + C_1\del J_0(\lam) + C_2\del^3 J_2(\lam) + C_3\del^5 J_4(\lam) + O(\del^7) \eeq
where again $\lam = b/\del$ and
 \beq J_n(\lambda) = \int_0^\infty
\frac{s_2(\sqrt{\sig^2+\lam^2})-1}{(\sig^2+\lam^2)^{3/2}}\sig^{n+3}\,d\sig 
   = \int_{|\lam|}^\infty  \frac{s_2(\rho)-1}{\rho^2} \sig^{n+2}\,d\rho \,,\quad \sig^2 = \rho^2 - \lam^2 \eeq
Proceeding as before, we seek a modified version of $s_2$ with the form
\beq s_2^{(7)}(\rho) = s_2(\rho) + b_1 \rho s_2'(\rho) + b_2 \rho^2 s_2''(\rho) + b_3 \rho^3 s_2'''(\rho) \eeq
To eliminate the error we need the $b$'s to satisfy the system
\beq \label{dblsys} b_1 J_{1n} + b_2 J_{2n} + b_3 J_{3n} = - J_n \,,\quad n=0,2,4 \eeq
with 
\beq J_{kn} = \int_{|\lam|}^\infty  \frac{\rho^ks_2^{[k]}(\rho)} {\rho^2} \sig^{n+2}\,d\rho \,,
\quad   \sig^2 = \rho^2 - \lam^2       \eeq
This is closely analogous to the single layer case.  (The extra minus sign in~\eqref{dblsys} is an artifact of our notation.)
Moreover, we find that
\beq -J_n = (n+2)I_n\,,\quad J_{kn} = (n+2)I_{kn}\,, \quad  n=0,2,4 \eeq
Thus the new system of equations is equivalent to \eqref{sys1}, and the solution agrees with the earlier one,
i.e., $b_k = a_k$ for $k=1,2,3$ with the $a$'s as in \eqref{tri1}--\eqref{tri3}.  After finding the derivatives
$s_2^{[k]}$ we obtain an expression for the seventh order version of $s_2$,
\beq 
	\label{s27} 
	s_2^{(7)}(\rho) = \erf(\rho) + \frac{2}{\sqrt{\pi}}  \left(-\rho +
  2(a_1 + 2a_2 + 2a_3)\rho^3 - 4(a_2 + 5a_3)\rho^5 + 8a_3\rho^7  \right)e^{-\rho^2}  
\eeq 
As before this formula reduces to the fifth order version, with $a_1^{(5)}$, $a_2^{(5)}$ as in \eqref{a5},
\beq  
	\label{s25}
	s_2^{(5)}(\rho) = \erf(\rho) + \frac{2}{\sqrt{\pi}}  \left(-\rho +
  2(a_1^{(5)} + 2a_2^{(5)})\rho^3  - 4a_2^{(5)} \rho^5 \right)e^{-\rho^2}  
\eeq 
as well as the third order version
 \beq \label{s23}  s_2^{(3)}(\rho) = \erf(\rho) + \left (- \frac{2}{\sqrt{\pi}}\rho  
  +  4I_0(\lam) e^{\lam^2}\rho^3 \right)e^{-\rho^2}  \eeq 
  
\medskip

{\bf The Stokeslet integral.} 
The equations of Stokes flow model incompressible fluid flow dominated by viscosity.  The velocity $\uu$ and pressure $p$
satisfy (e.g. see \cite{pozbook})
\beq \label{stokeseqns}
  - \Delta\uu + \nabla p = 0\,,\quad \nabla\cdot\uu = 0 \eeq
in the absence of force with viscosity set to $1$.
The fundamental solution for the velocity is the Stokeslet
\beq
	\label{Stokeslet}
	S_{ij}(\bd{y,x}) = \frac{\del_{ij}}{|\bd{y} - \bd{x}|} + \frac{(y_i - x_i
)(y_j - x_j)}{|\bd{y} - \bd{x}|^3}
\eeq	
with $i,j = 1,2,3$.  A force $\ff$ on a surface
$\Gamma$ determines the velocity
\beq \label{stosgl}
	u_i(\bd{y}) = \frac{1}{8\pi}\int_\Gamma S_{ij}(\bd{y,x}) f_j(\bd{x})
dS(\bd{x}) \eeq
A sum over $j$ is implicit on the right side.  Again we use a subtraction.
If $f_j$ in \eqref{stosgl} is the normal vector $n_j$, then the integral is zero
(see \cite{pozbook}, eqns. (2.1.4) and (6.4.3)).
With $\xx_0$ as before we can rewrite \eqref{stosgl} as
\beq \label{stosglsub}
u_i(\bd{y}) = \frac{1}{8\pi}\int_\Gamma S_{ij}(\bd{y,x})
 [f_j(\xx) - f_k(\xx_0)n_k(\xx_0)n_j(\xx)]\,dS(\bd{x}) \eeq
To compute \eqref{stosglsub} we replace $S_{ij}$ with the regularized version
\beq \label{Sreg} S_{ij}^\del(\bd{y,x}) = \frac{\del_{ij}}{r}s_1(r/\del)
     + \frac{(y_i - x_i)(y_j - x_j)}{r^3}s_2(r/\del)\,,
         \quad r = |\yy - \xx| \eeq
with $s_1$ and $s_2$ as in \eqref{sglshape},\eqref{dblshape}, resulting in a smooth kernel.
The higher order versions of $s_1$, $s_2$ lead to high order versions of \eqref{stosglsub}.

\medskip

{\bf The stresslet integral.}
The stresslet integral or double layer integral in Stokes flow is
\beq \label{stodbl}
v_i(\bd{y}) = \frac{1}{8\pi}\int_\Gamma T_{ijk} (\bd{y,x}) q_j(\bd{x}) n_k(\bd{x})dS(\bd{x}) \eeq
with kernel
\beq \label{stresslet}
T_{ijk} (\bd{y,x}) = -\frac{6(y_i - x_i)(y_j - x_j)(y_k - x_k)}{|\bd{y} -\bd{x}|^5} \eeq
We use an identity (see \cite{pozbook}, eqns. (2.1.12) and (6.4.5)) to write \eqref{stodbl}
in the subtracted form, with $\chi$ as before,
\beq \label{stodblsub}
	v_i(\bd{y}) = \frac{1}{8\pi}\int_\Gamma T_{ijk} (\bd{y,x}) [q_j(\bd{x})
	 - q_j(\bd{x}_0)] n_k(\bd{x}) dS(\bd{x}) + \chi (\bd{y}) q_i(\bd{x}_0)  \eeq
We need to rewrite the kernel \eqref{stresslet} to make it
accessible to our method.
For $\yy$ near the surface we have, as in \eqref{closest},
$\yy = \xx_0 + b\nn$ with $\xx_0 \in \Gamma$ and $\nn = \nn(\xx_0)$.
In the numerator of $T_{ijk}$ we substitute $\yy - \xx = b\nn - \xxh$, with $\xxh = \xx - \xx_0$.
In terms with factors of $b^2$ we substitute $b^2/r^2 = 1 - (r^2 - b^2)/r^2$ where
$r = |\yy - \xx|$ and then $r^2 - b^2 = |\xxh|^2 - 2b\xxh\cdot\nn$.  In this way we split
the kernel as $T_{ijk} =  T_{ijk}^{(1)} + T_{ijk}^{(2)}$ so that $T_{ijk}^{(1)}$ has denominator
$r^3$ and $T_{ijk}^{(2)}$ has denominator $r^5$.  Specific formulas are given in \cite{extrap}, eqns. (23)--(26).

To compute \eqref{stodblsub} we replace $T_{ijk}$ with the regularized version 
\beq \label{Treg} 
T_{ijk}^\del = T_{ijk}^{(1)}s_2(r/\del) + T_{ijk}^{(2)}s_3(r/\del) \eeq
where
\beq \label{s3}
     s_3(\rho) = \erf(\rho) 
            - \frac{2}{\sqrt{\pi}}\left(\frac23 \rho^3 + \rho\right)e^{-\rho^2}  \eeq
For the first part of \eqref{Treg} we already have high order versions of $s_2$, whereas $s_3$ in the
second part is new.  In \cite{extrap} we found that the error in the integral with 
$T_{ijk}^{(2)}s_3$ is of the form
\beq C_1\del K_0(\lam) + C_2\del^3 K_2(\lam) + C_3\del^5 K_4(\lam) + O(\del^7) \eeq
with
\beq K_n(\lambda) = \int_0^\infty
\frac{1 - s_3(\sqrt{s^2+\lam^2})}{(\sig^2+\lam^2)^{5/2}}\sig^{n+5}\,d\sig 
 = \int_{|\lam|}^\infty 
 \frac{1 - s_3(\rho)}{\rho^4}\sig^{n+4}\,d\rho \,,\quad \sig^2 = \rho^2 - \lam^2   \eeq
As before we look for a seventh order version $s_3^{(7)}$ of $s_3$ in the form
\beq 
   s_3^{(7)}(\rho) = s_3(\rho) + c_1\rho s_3'(\rho) + c_2\rho^2 s_3''(\rho) +c_3\rho^3 s_3'''(\rho) \eeq
The condition to remove the leading error is
\beq c_1 K_{1n} + c_2 K_{2n} + c_3 K_{3n} = K_n\,,\quad  n = 0,2,4 \eeq
where
\beq K_{jn}(\lam) = \int_{|\lam|}^\infty 
 \frac{\rho^j s_3^{[j]}(\rho)}{\rho^4}\sig^{n+4}\,d\rho \,,\quad \sig^2 = \rho^2 - \lam^2   \eeq
 Remarkably, we find that
\beq K_n = \kappa_n I_n\,,\quad K_{jn} = \kappa_n I_{jn}\,, \quad  n=0,2,4;\; j=1,2,3 \eeq
with $\kappa_0 = 8/3$  $\kappa_2 = 8$  $\kappa_4 = 16$,
so that once again the system of equations is equivalent to \eqref{sys1}, and the coefficients are given by 
\eqref{tri1}--\eqref{tri3}.  After calculating derivatives of $s_3$ we find the high order
versions of $s_3$,
\beq 
	\label{s37} 
	s_3^{(7)}(\rho) = s_3(\rho) + \frac{8}{3\sqrt{\pi}} \left(
   (a_1 + 4a_2 + 12a_3)\rho^5 - 2(a_2 + 9a_3)\rho^7 + 4a_3\rho^9 \right)e^{-\rho^2}  
\eeq
\beq 
	\label{s35} 
	s_3^{(5)}(\rho) = s_3(\rho) + \frac{8}{3\sqrt{\pi}} \left(
(a_1^{(5)} + 4a_2^{(5)} )\rho^5 - 2a_2^{(5)}\rho^7 \right)e^{-\rho^2}  
\eeq
with $a_1^{(5)}$, $a_2^{(5)}$ as in \eqref{a5}, and the third order version
\beq 
	\label{s33} 
	s_3^{(3)}(\rho) = \erf(\rho) + \frac{2}{\sqrt{\pi}}\left(- \rho  - \frac23\rho^3  
      + \frac43\sqrt{\pi} I_0(\lam)e^{\lam^2}\rho^5 \right)e^{-\rho^2}  
\eeq
      
\medskip

{\bf Pressure in Stokes flow.}
The pressure due to a force on a surface $\Gamma$ in Stokes flow is
\beq \label{pressure} p(\yy) = 
          \int_\Gamma \nabla G(\yy-\xx)\cdot \ff(\xx)\,dS(\xx)  \eeq 
with $G$ as in \eqref{sgllayer}.
We can separate the force $\ff(\xx)$ into normal and tangential parts, 
\beq \ff = f^{(n)}\nn + \ff^{(t)} \equiv (\ff\cdot\nn)\nn - \nn\times(\nn\times\ff) \eeq
where $\nn = \nn(\xx)$ is the outward normal. The pressure is now a sum of two parts, $p = p^{(n)} + p^{(t)}$.
The first part amounts to a double layer potential as already discussed.
We can subtract as before, with $\xx_0$ as in \eqref{closest},
\beq \label{pn} p^{(n)}(\yy) = - \int_\Gamma \frac{\pa G(\xx-\yy)}{\pa\nn(\xx)}
                   [f^{(n)}(\xx) - f^{(n)}(\xx_0)]\,dS(\xx) - f^{(n)}(\xx_0)\chi(\yy) \eeq
To compute $p^{(n)}$ we regularize with a version of $s_2$.
For the tangential part the integrand is
\beq \nabla G \cdot [- \nn\times (\nn\times \ff)] =
   (\nn\times \nabla G)\cdot(\nn\times \ff) \eeq
We can make a similar subtraction in this part using the identity (e.g., see \cite{pozjcp}, page 284) 
\beq \int_\Gamma \nn(\xx)\times \nabla G(\yy-\xx) \,dS(\xx) = 0 \eeq
We can use this fact to write
\beq p^{(t)}(\yy) = \int \left(\nn(\xx)\times \nabla G(\yy-\xx)\right) \cdot 
         \left[\nn(\xx)\times \ff(\xx) - \nn(\xx_0)\times \ff(\xx_0)\right] \,dS(\xx) \eeq
and we can again regularize with a version of $s_2$ as previously described.


\section{Formulas for points on the surface}
\label{sec:onSurf}

We summarize formulas for the important special case of evaluation on the surface.  All the formulas in the
previous section are valid on the surface, but some simplifications are possible in this case.
For $\yy$ on the surface $\lam = 0$, and the equations \eqref{tri1}--\eqref{tri3} reduce to
\beq a_3 = \frac{1}{15}\,;\quad a_2 = \frac13 + 7a_3\,; \quad a_1 = 1 + 2a_2 - 6a_3 \eeq
We determine the $a$'s from these equations for seventh order and those for the lower orders
as before. With order $p = 7$, $5$, $3$ we get
\beq \label{aon7} a_1 = 11/5\,,\quad a_2 = 4/5\,,\quad a_3 = 1/15\,,\qquad p = 7 \eeq
\beq \label{aon5} a_1 = 5/3\,,\quad a_2 = 1/3\,,\quad a_3 = 0\,,\qquad p = 5 \eeq
\beq \label{aon3} a_1 = 1\,,\quad a_2 = a_3 = 0\,,\qquad p=3 \eeq

\medskip

{\bf The single layer potential.}  For the single layer potential \eqref{sglreg} evaluated on the surface, we
combine \eqref{aon7}--\eqref{aon3} with \eqref{s17}--\eqref{s13}.  We get 
\beq s_1^{(p)}(\rho) = \erf(\rho)  + \frac{2}{\sqrt{\pi}} e^{-\rho^2} m(\rho) \eeq
with
\beq \label{s1_onSurf_7} m(\rho) = \frac{11}{5}\rho - \frac{26}{15}\rho^3 + \frac{4}{15}\rho^5\,,\quad p = 7 \eeq
\beq \label{s1_onSurf_5} m(\rho) = \frac53\rho - \frac23\rho^3\,,\quad p = 5\,; 
           \qquad  m(\rho) = \rho\,,\quad p = 3 \eeq
The fifth order formula was used in \cite{byw}, eqn. (3.14) and \cite{tbjcp}, eqn. (50).

In the single layer or Stokeslet integral we have $s_1(\rho)/r = \delta^{-1}s_1(\rho)/\rho$, and the value is needed at $\rho = 0$.  As $\rho \to 0$, $\erf(\rho)/\rho \to 2/\sqrt{\pi}$, and for $p = 3$, $5$, or $7$, $s_1^{(p)}(\rho)/\rho \to (2/\sqrt{\pi})(1 + a_1)$, with $a_1$ as in \eqref{aon7}--\eqref{aon3}.

\medskip

{\bf The double layer potential.} The expressions \eqref{s27}--\eqref{s23}
already derived for the double layer can be used on the surface,
but a simplification can be made in this special case with the integral in the subtracted form
\eqref{dblsub}.  On the surface there is no contribution to the error \eqref{dblerr} from $J_0$, and
$s_2$ as defined in \eqref{dblshape} is third order.
To find the seventh order version, we can drop the last term in \eqref{dblsys} and the equation with $n=0$.
The system of equations \eqref{dblsys} for the coefficients then reduces to 
\beq  b_1 I_{1n} + b_2 I_{2n} = I_n \,,\quad n=2,4 \eeq
We find that $b_1 = 3/5$, $b_2 = 1/15$, and in place of \eqref{s27} we have
\beq s_{2*}^{(7)} = \erf(\rho)  + \frac{2}{\sqrt{\pi}} e^{-\rho^2} 
  \left[-\rho + \frac{22}{15}\rho^3 - \frac{4}{15}\rho^5 \right] \eeq
Similarly the fifth order version is, with $b_1 = 1/3$, $b_2 = 0$,
\beq s_{2*}^{(5)} = \erf(\rho)  + \frac{2}{\sqrt{\pi}} e^{-\rho^2} 
  \left[-\rho + \frac{2}{3}\rho^3 \right] \eeq
These formulas can be used with \eqref{dblreg} to evaluate the double layer potential on the surface, with $\chi = 1/2$.  
The fifth order version was used in \cite{b04}, eqn. (1.12) and in \cite{byw}, eqn. (3.17).

\medskip

{\bf The Stokeslet integral.}  We regularize the Stokeslet as in \eqref{Sreg} with versions of $s_1$ and $s_2$.
Versions of $s_1$ on the surface were discussed above.  For $s_2$ we can proceed similarly.  
(We cannot use the special form of $s_2$ described above for the double layer.)
We combine 
\eqref{aon7}--\eqref{aon3} with \eqref{s27}--\eqref{s23} and obtain
\beq s_2^{(p)}(\rho) = \erf(\rho)  + \frac{2}{\sqrt{\pi}} e^{-\rho^2} m(\rho) \eeq
with
\beq 
	\label{s2_onSurf_7} 
	m(\rho) = - \rho + \frac{118}{15}\rho^3 - \frac{68}{15}\rho^5 + \frac{8}{15}\rho^7\,,\quad p = 7 
\eeq
\beq 
	\label{s2_onSurf_5} 
	m(\rho) = - \rho + \frac{14}{3}\rho^3 - \frac43\rho^5\,,\quad p = 5\,; 
           \qquad  m(\rho) = - \rho + 2\rho^3\,,\quad p = 3 
\eeq
For evaluation on the surface the subtraction in
\eqref{stosglsub} is not necessary; the regularization can be used with or without it.
That is, we can regularize \eqref{stosgl}, rather than \eqref{stosglsub}, using \eqref{Sreg}.
The fifth order version of $s_2$ was used in \cite{tbjcp}, eqn. (51), and the fifth order versions
of $s_1$ and $s_2$ were used in \cite{agbiros}.

\medskip

{\bf The stresslet integral.}  To evaluate the stresslet integral on the surface in the subtracted form \eqref{stodblsub}
we do not need the splitting \eqref{Treg} of the kernel $T_{ijk}$ described in Sect. 2.  We only need to insert a version of
$s_3$ in \eqref{stodblsub}, with $\chi = 1/2$.    As noted for the harmonic double layer, there is no contribution
to the error from the $K_0$ term.  Thus we can simplify as we did in the double layer case.  With the coefficients
$b_1$, $b_2$ as above we obtain from \eqref{s37} and \eqref{s35} the seventh and fifth order versions
\beq 
	\label{s3_onSurf_7}
	s_{3*}^{(7)} = \erf(\rho)  + \frac{2}{\sqrt{\pi}} e^{-\rho^2} 
  \left[-\rho - \frac23\rho^3 + \frac{52}{45}\rho^5  - \frac{8}{45}\rho^7 \right] 
\eeq
\beq 
	\label{s3_onSurf_5}
	s_{3*}^{(5)} = \erf(\rho)  + \frac{2}{\sqrt{\pi}} e^{-\rho^2} 
  \left[-\rho - \frac23\rho^3 + \frac49\rho^5\right] 
\eeq
The fifth order version was used in \cite{novel}, eqn. (14).

\bigskip


\section{Numerical examples}
We present examples to illustrate the use of the regularization formulas.  We imbed a surface in a cubic box and discretize with a regular grid of size $h$.  Regularized integrals are computed with sums over quadrature points where the surface intersects grid intervals; further description of the quadrature is given in \cite{extrap}, Sect.4 or \cite{byw}.  We compute the integrals at grid points within distance $O(h)$ of the surface.  We neglect the regularization for quadrature points at distance $8\del$ or more from the target point, since its effect would be negligible.  We report maximum and $L^2$ errors, with the 
$L^2$ norm defined as usual,
\beq \|\epsilon\|_{L^2} = \left( \sum |\epsilon(\yy)|^2 / N \right)^{1/2} \eeq
where $\epsilon(\yy)$ is the error at $\yy$ and $N$ is the number of points. 
We present absolute errors.  

In Sect.~\ref{sec:test1_lapl} we verify the accuracy of single and double layer potentials on a sphere.  Examples in Sect.~\ref{sec:test2} with $\delta = \kappa h^q$ lead to a strategy for choosing the parameter $\kappa$.  In Sect.~\ref{sec:lapl_grid} we describe the procedure for obtaining values of a harmonic function at all grid points from integral values at points near the interface by forming a discrete fourth order Laplacian and inverting, as in~\cite{mayo85}.  In Sect.~\ref{sec:stokes_grid} we extend this procedure to solve for the pressure and velocity in Stokes flow around a translating spheroid.  Finally in Sec.~\ref{sec:stokes_2surf} we solve an integral equation for the velocity in Stokes flow past two spheroids with surface tension which are close together.  In all the examples but the last, computational results are compared with known solutions.  For the two spheroids we use a treecode~\cite{wang-krasny-tlupova-20} for efficient summation with parameters determined in~\cite{siebor-tlupova}. 


\subsection{Single and double layer potentials on a sphere and the choice of $\delta$ }  
\label{sec:test1_lapl}

We begin with solutions of Laplace's equation determined by a spherical harmonic
on the unit sphere centered at the origin.
We represent the solution by either a single or double layer potential,
\eqref{sgllayer} or \eqref{dbllayer}.  For a general
surface $\Gamma$ the two potentials are characterized by the jump conditions 
\beq \label{jump1} [\Sl(\xx)] = 0\,,\quad [\pa \Sl(\xx)/\pa\nn] = f(\xx) \eeq
\beq \label{jump2} [\Dl(\xx)] = -g(\xx)\,,\quad [\pa \Dl(\xx)/\pa\nn] = 0 \eeq
where $[\cdot]$ means the value outside $\Gamma$ minus the value inside.

On the unit sphere we choose the spherical harmonic function
\beq f(\xx) = (7/8)(x_1 - 2x_2)(15x_3^2 - 3)\,, \quad |\xx| = 1 \eeq
The functions
\beq u_-(\yy) = r^3f(\yy/r)\,,\quad u_+(\yy) = r^{-4}f(\yy/r)\,,
    \quad\; r = |\yy| \eeq
are harmonic inside and outside, respectively.  We define
$\Sl(\yy)$ by \eqref{sgllayer} and $\Dl(\yy)$ by \eqref{dbllayer}
with $g = f$.  The jump conditions determine them as
\beq \Sl(\yy) = - (1/7)u_-(\yy)\,, \quad |\yy|<1\,; \quad\;
   \Sl(\yy) = - (1/7)u_+(\yy)\,, \quad |\yy|>1\eeq
\beq \Dl(\yy) =  (4/7)u_-(\yy)\,, \quad |\yy|<1\,; \quad\;
   \Dl(\yy) = - (3/7)u_+(\yy)\,, \quad |\yy|>1\eeq
   
\begin{figure}[!htb]
\centering
\scalebox{0.425}{\includegraphics{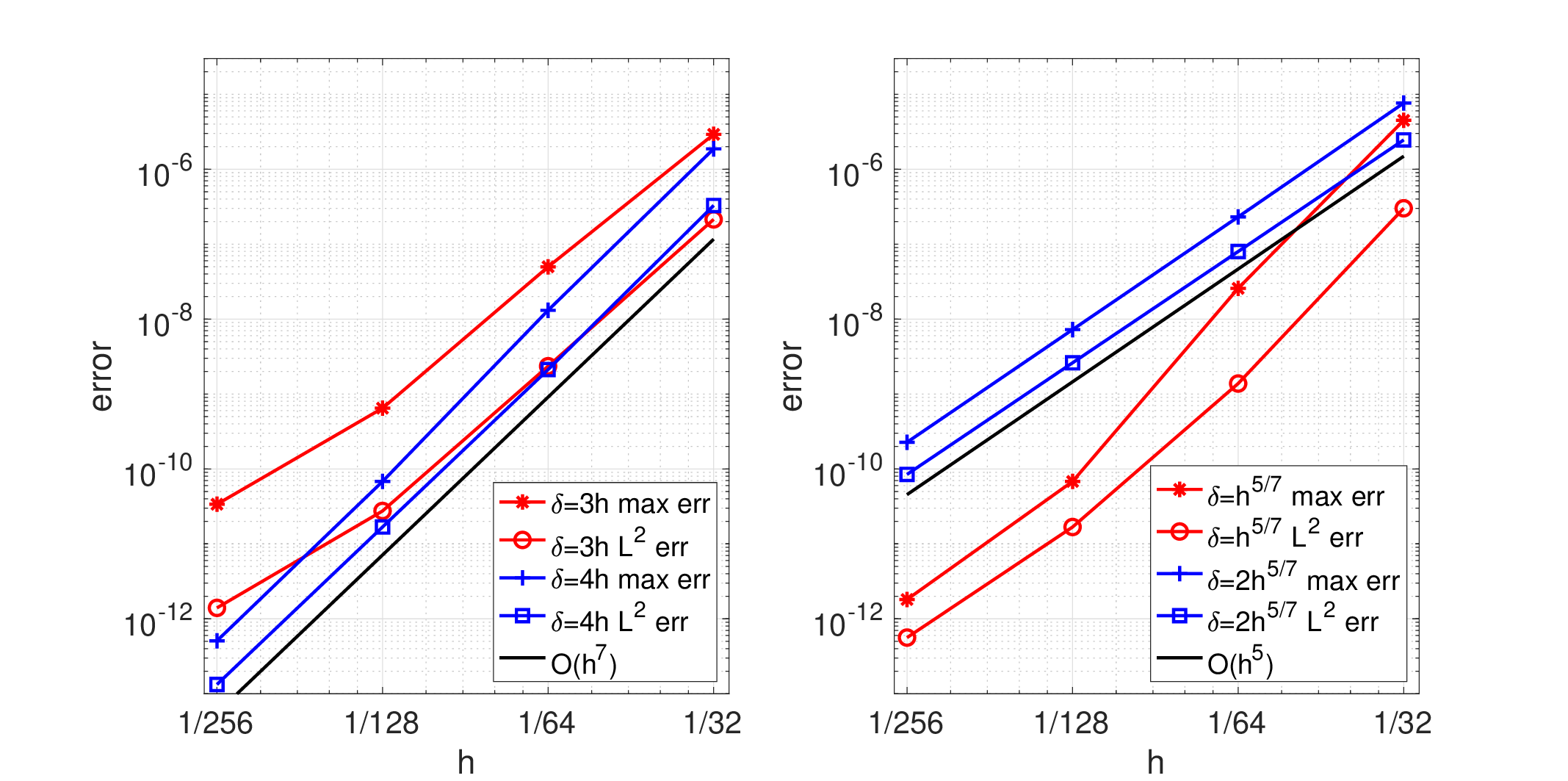}} 
\caption{Errors for the harmonic single layer on a unit sphere, at grid points within distance $h$ to
 the sphere.}
\label{figure:lapl_SL_sphere}
\end{figure}

\begin{figure}[!htb]
\centering
\scalebox{0.425}{\includegraphics{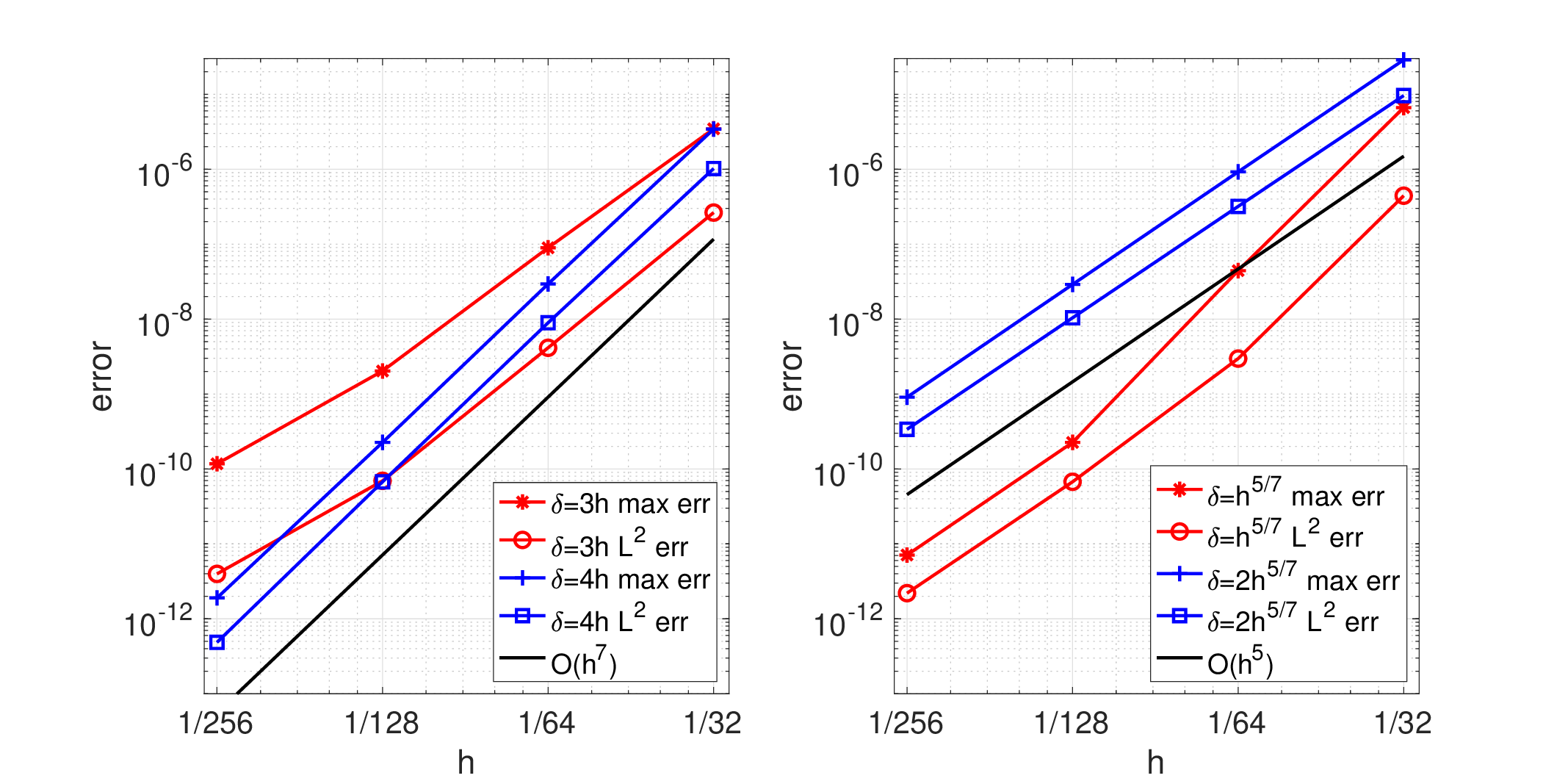}} 
\caption{Errors for the harmonic double layer on a unit sphere, at grid points within distance $h$ to
 the sphere.}
\label{figure:lapl_DL_sphere}
\end{figure}

The known solutions enable us to test the accuracy of the single and double layer potentials separately,
and especially the dependence of the error on the relation between $\delta$ and $h$.  We use
the seventh order versions of $s_1$ and $s_2$ given in~\eqref{s17} and~\eqref{s27}. We have chosen $\delta = 3h$, $\delta = 4h$, $\delta = h^{5/7}$ and $\delta = 2h^{5/7}$. We report the maximum and $L^2$ errors for the integrals at grid points within distance $h$ of the sphere.  (Any point that happened to be on the sphere was excluded.)  
The maximum and $L^2$ norms for $|\Sl|$ were about $1.15$ and $.50$.  For $|\Dl|$ they were $4.6$ and $1.8$.
Figure~\ref{figure:lapl_SL_sphere} displays results for the single layer potential, and Figure~\ref{figure:lapl_DL_sphere} has those for the double layer potential.  The results in the two cases are similar.  With $\delta = 4h$ we see $O(h^7)$ error, as expected; for $\delta = 3h$ the predicted order holds for larger $h$ but fails for the smallest $h$.
In contrast, with $\delta = 2h^{5/7}$ we see the expected error $O(h^5)$, while for
$\delta = h^{5/7}$ the error shows the expected error for small $h$ and decreases more rapidly for the larger values of $h$.  This behavior is consistent with our expectation that making $\delta/h$ constant can give the expected order for moderate $h$, but not in the limit $h \to 0$ because the discretization error becomes significant; with $\delta$ proportional to $h^{5/7}$ it appears the smoothing error remains dominant as $h \to 0$.  
  
\medskip


\subsection{Convergence with $\delta = \kappa h^q$}
\label{sec:test2}

We test the accuracy of the various regularized kernels of order $p$,
choosing $\del = \kap h^q$.  The predicted order of accuracy is $O(h^{pq})$.  For $p=3$, we choose
$q = 2/3$; for $p=5$, $q = 4/5$; and for $p=7$, $q = 5/7$.  
To choose $\kap$, we first choose $\del$ for $h = 1/64$, a relatively coarse value of $h$.  We define
$\kapz = \del/h$ for this case.  Then for all $h$, we set $\del = \kap h^q$ with $\kap$ determined
so that $\del/h = \kapz$ at $h = 1/64$; that is, $\kap = \kapz/64^{1-q}$.
The test problems are known solutions.  We evaluate the integrals at a random selection of
grid points within $O(h)$ of the surface; the random choice allows us to observe the trend for
small $h$.

\begin{figure}[!htb]
\centering
\scalebox{0.64}{\includegraphics{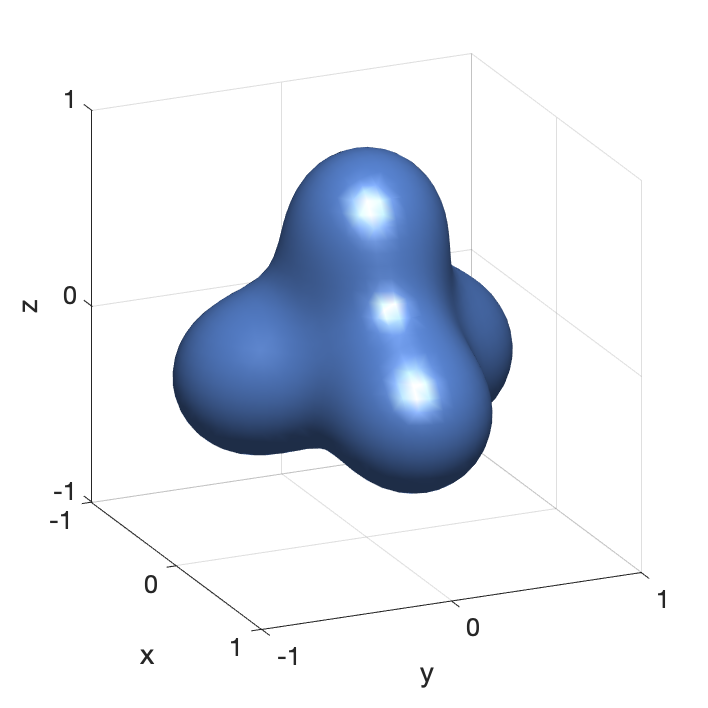}} 
\scalebox{0.375}{\includegraphics{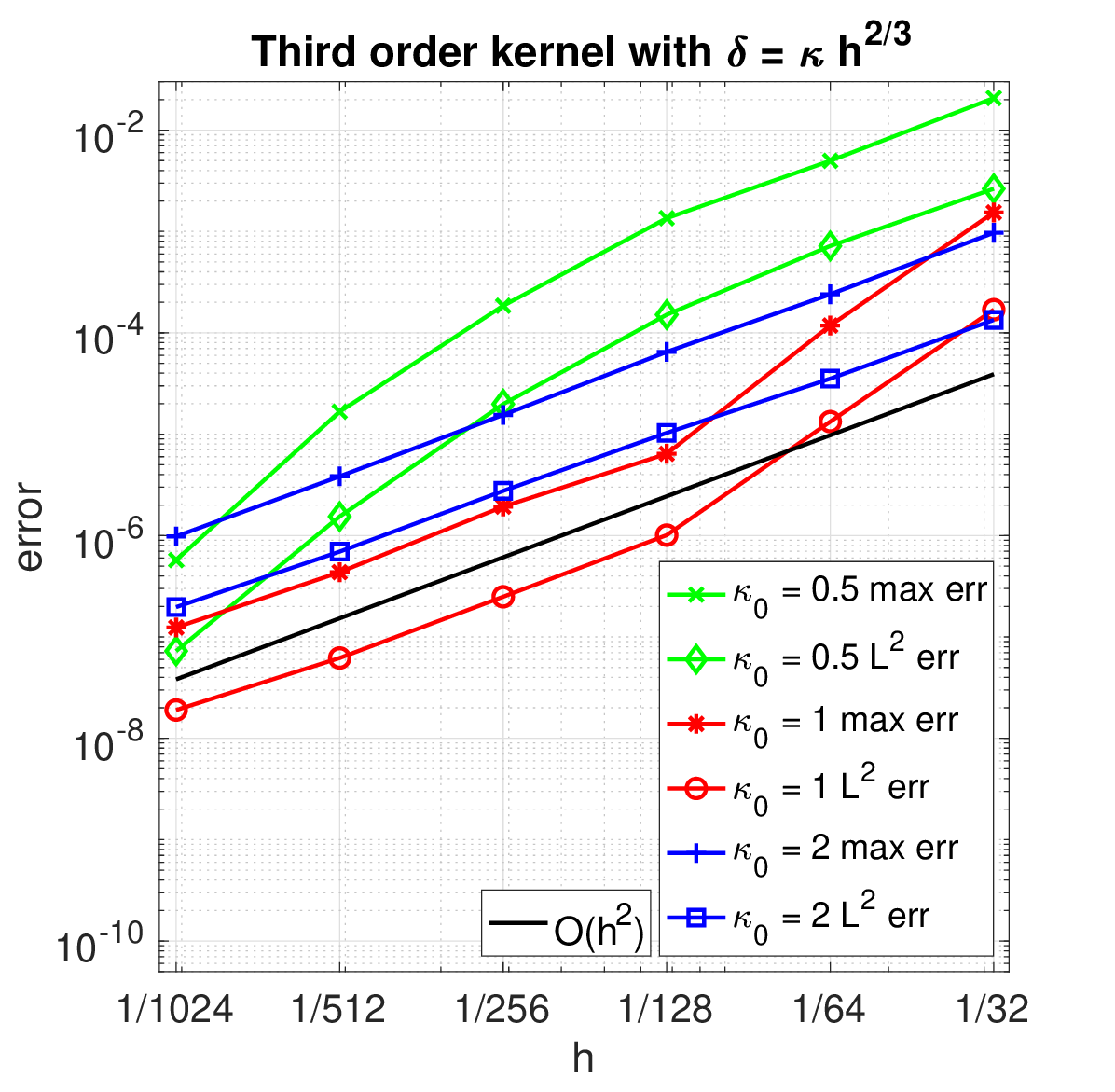}} 
\scalebox{0.375}{\includegraphics{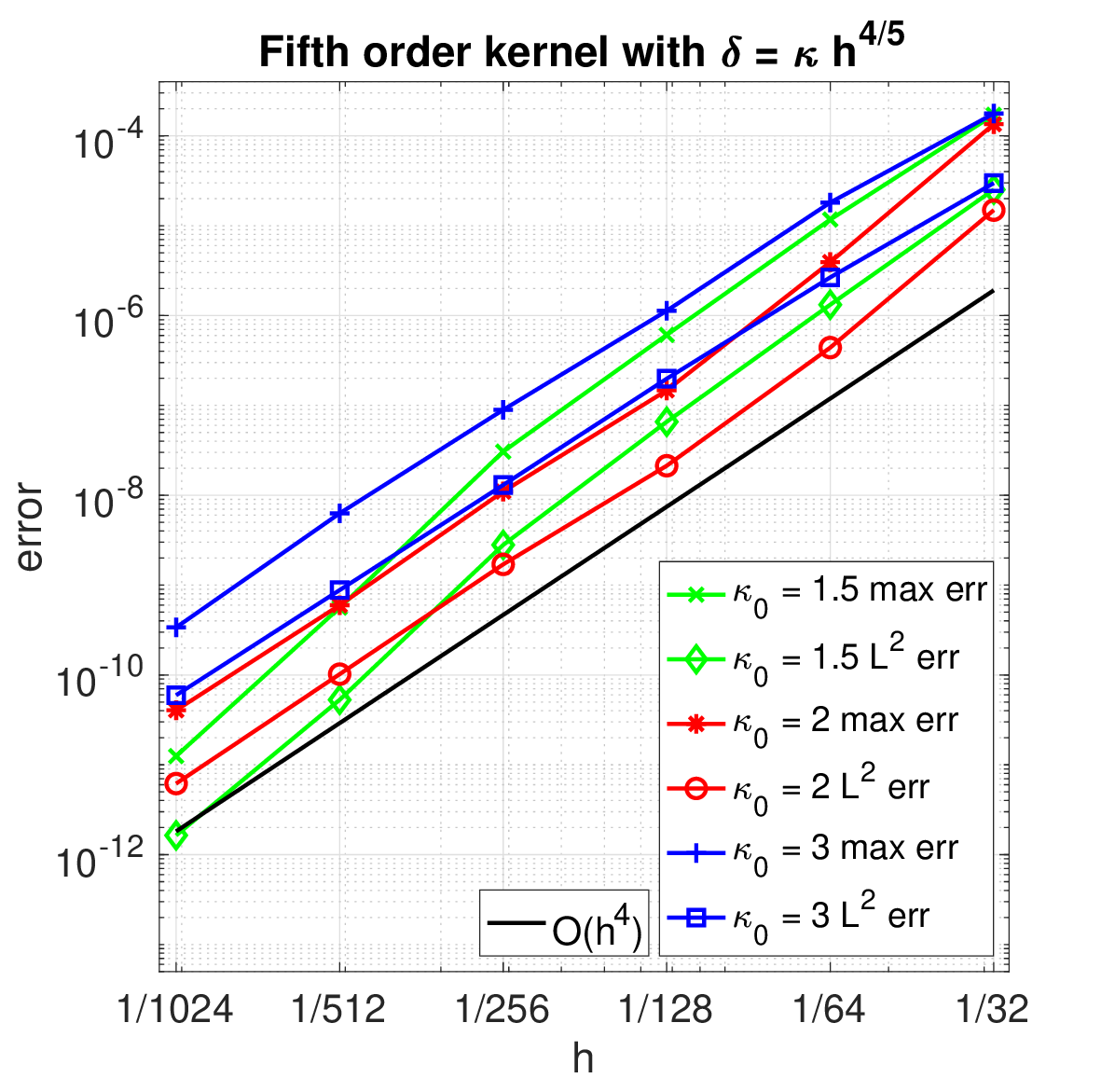}} 
\scalebox{0.375}{\includegraphics{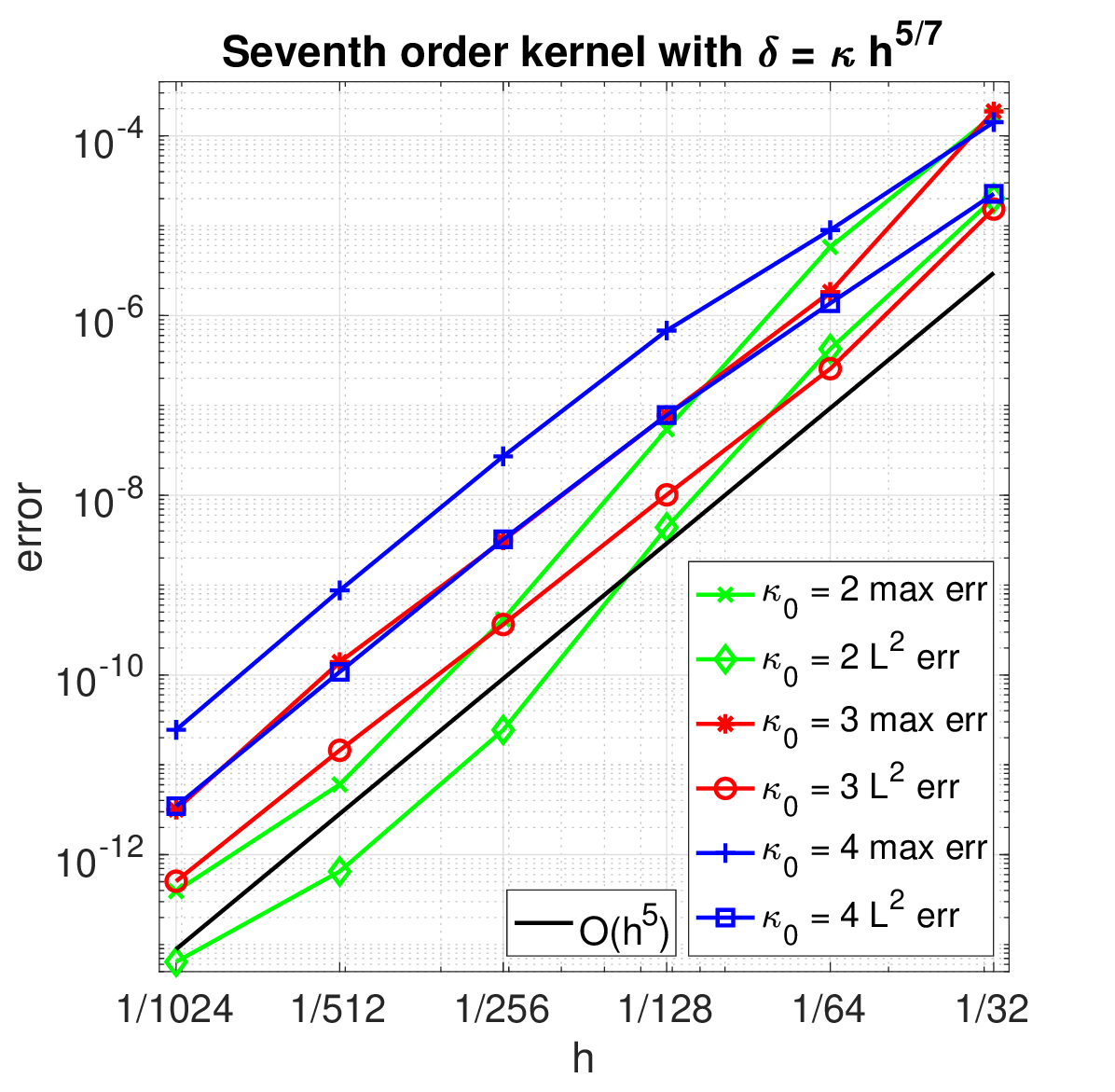}} 
\caption{Errors for the harmonic solution on the molecular surface, at grid points randomly chosen within distance $h$ to the surface. Here $\kappa$ is chosen so that $\delta/h = \kappa_0$ when $h = 1/64$.}
\label{figure:lapl_mol}
\end{figure}

For our first tests we choose harmonic functions $u_+$ outside and $u_-$ inside.  We set $f = [\pa u/\pa n]$ and $g = -[u]$, the jumps across $\Gamma$ as above.  Then assuming $u_+$ decays at infinity, $ u(\yy) = \Sl(\yy) + \Dl(\yy)$ on both sides, in view of the jump conditions \eqref{jump1}, \eqref{jump2}, where $\Sl$ and $\Dl$ are defined in \eqref{sgllayer}, \eqref{dbllayer}.  We choose
\beq \label{basicu} u_-(\yy) = (\sin{y_1} + \sin{y_2})\exp{y_3}\,,
   \quad u_+(\yy) = 0  \eeq
Our surface represents a molecule with four atoms, pictured in Figure~\ref{figure:lapl_mol},
\beq \label{molesurf} \sum_{i=1}^4 \exp(-|\xx - \xx_k|^2/r^2) = c \eeq
with $r = .5$, $c = .6$, and $\xx_k$ given by
\beq (\sqrt{3}/3,0,-\sqrt{6}/12)\,,\; 
   (-\sqrt{3}/6,\pm .5,-\sqrt{6}/12)\,,\;
   (0,0,\sqrt{6}/4)  \eeq
The solution is computed as $\Sl(\yy) + \Dl(\yy)$, with regularization, at grid points within $h$ of the surface.  The points are selected randomly with probability $85h^2$, resulting in a choice of about $1200$ points in each case.  For each order $p = 3$, $5$, $7$, maximum and $L^2$ errors are reported in Figure~\ref{figure:lapl_mol}.
For each $p$ three choices of $\kapz$ are shown.  In each case the predicted order $O(h^{pq})$ can be seen
for the larger two values of $\kapz$ but is less clear for the smallest values.

\begin{figure}[!htb]
\centering
\scalebox{0.375}{\includegraphics{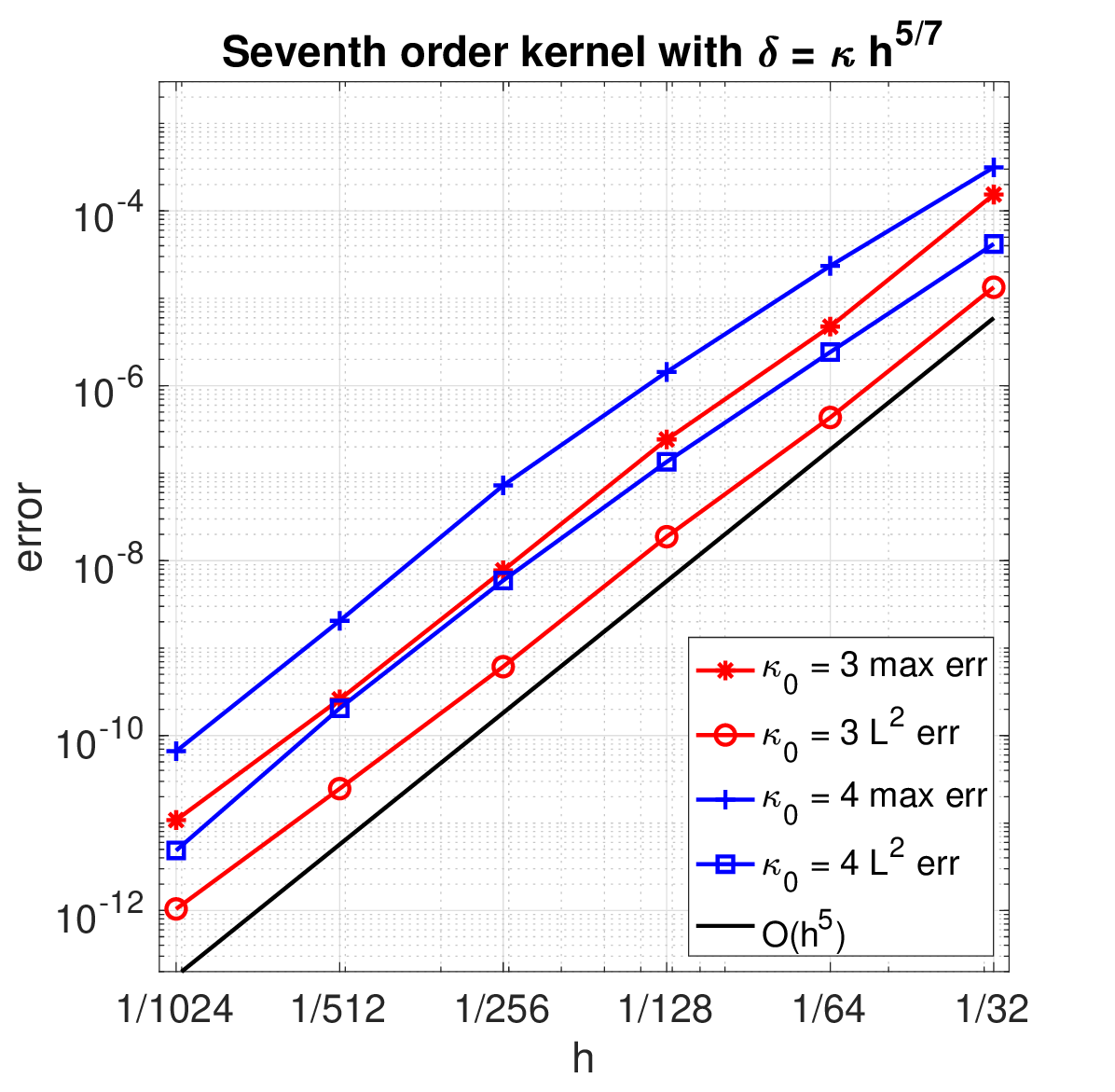}} 
\scalebox{0.375}{\includegraphics{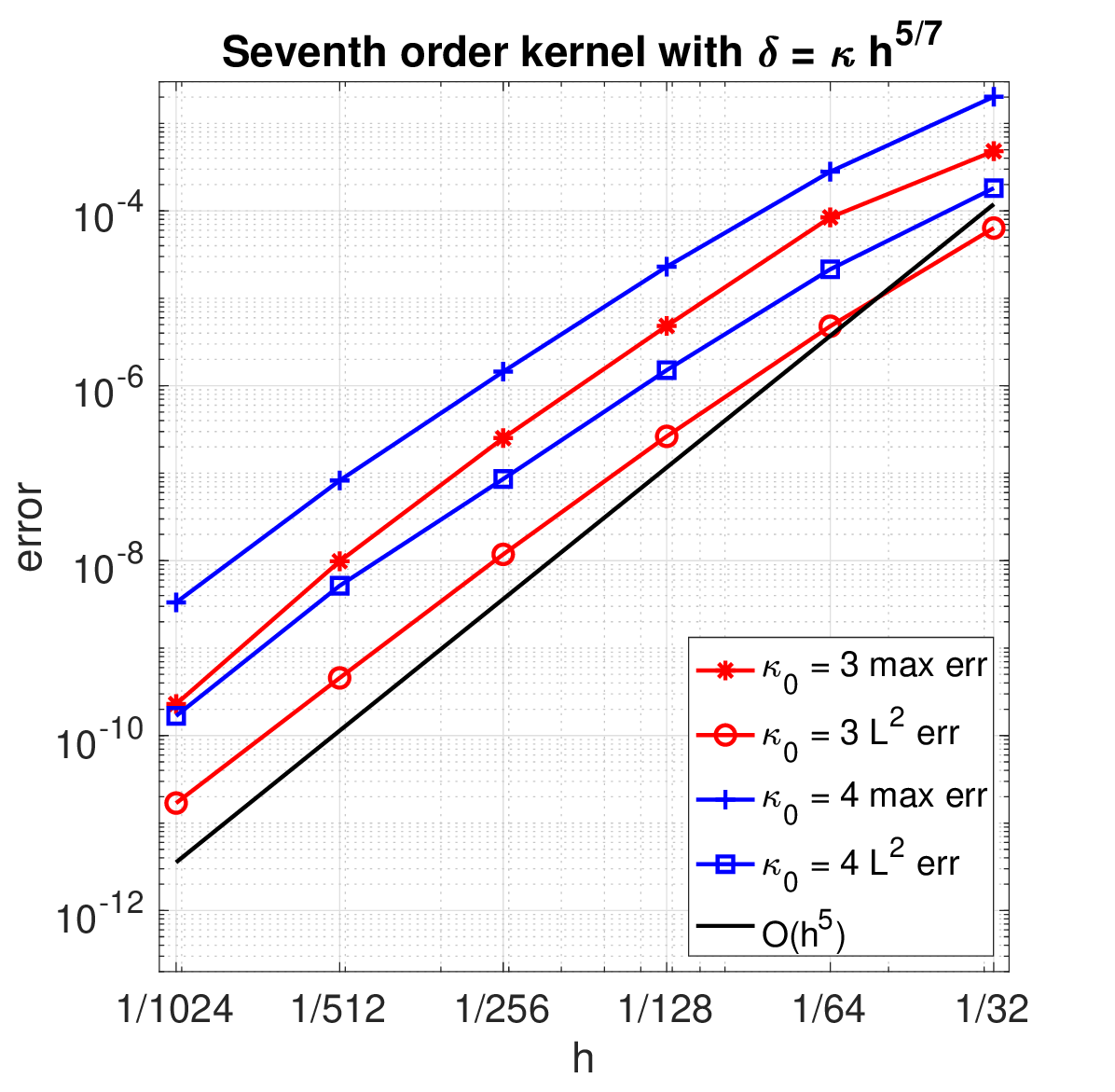}} 
\caption{Errors for the harmonic solution on an ellipsoid (1,\;.6,\;.4) (left) and ellipsoid (1,\;.4,\;.3) (right), at grid points randomly chosen within distance $h$ to the surface. Here $\kappa$ is chosen so that $\delta/h = \kappa_0$ when $h = 1/64$.} 
\label{figure:lapl_ellipsoid}
\end{figure}

We perform a similar test with two ellipsoids and $p = 7$, $q = 5/7$, $\kapz = 3$ or $4$.
The ellipsoids have semi-axes $1$, $.6$, $.4$ and $1$, $.4$, $.3$,
and the exact solution is again \eqref{basicu}.  We select target grid
points within $h$ of the ellipsoid, with probability $120h^2$ or $180h^2$, resp., amounting to a choice of
about $1300$ points.  Results are shown in
Figure~\ref{figure:lapl_ellipsoid}.  They show the expected $O(h^5)$ accuracy with $\kapz = 4$, somewhat less so for $\kapz = 3$.

\begin{figure}[!htb]
\centering
\scalebox{0.64}{\includegraphics{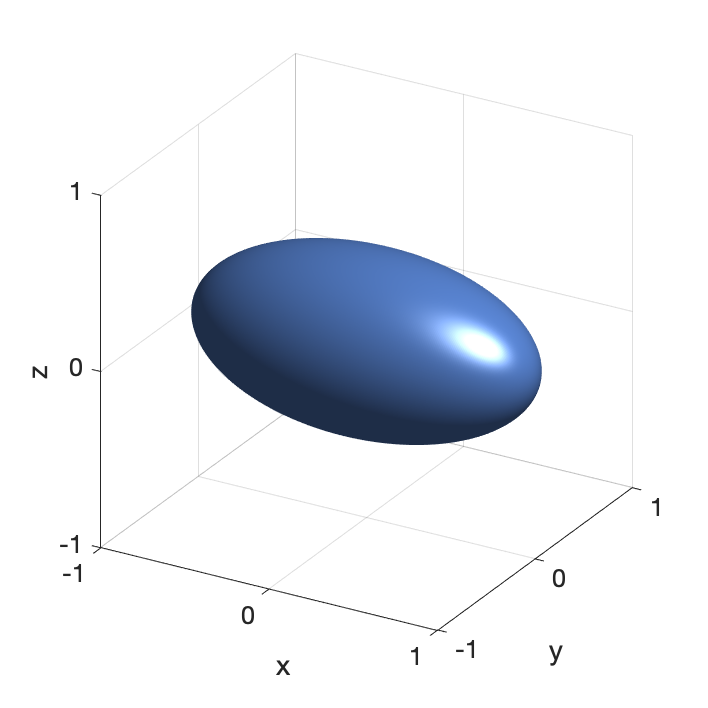}} 
\scalebox{0.375}{\includegraphics{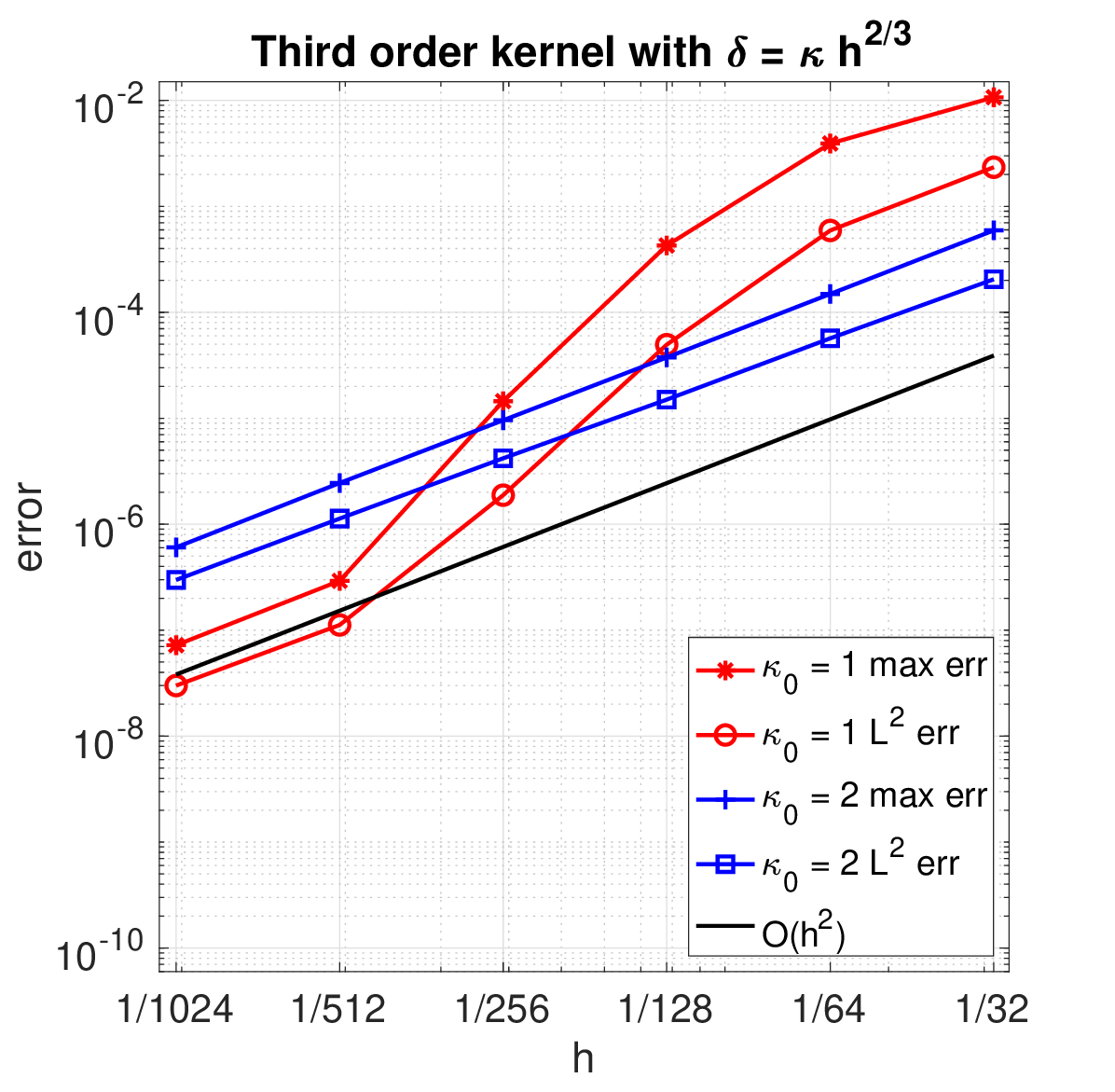}} 
\scalebox{0.375}{\includegraphics{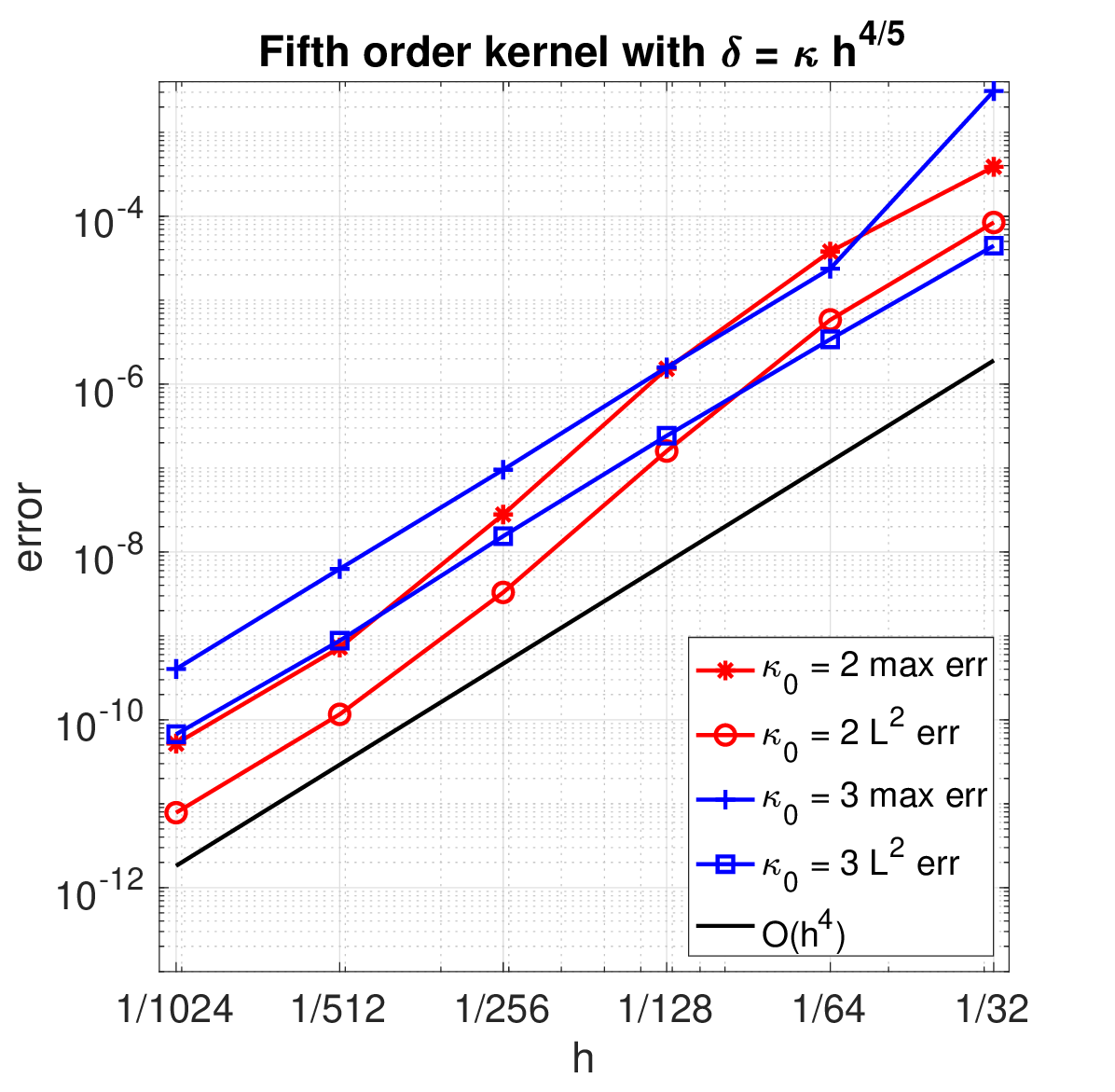}} 
\scalebox{0.375}{\includegraphics{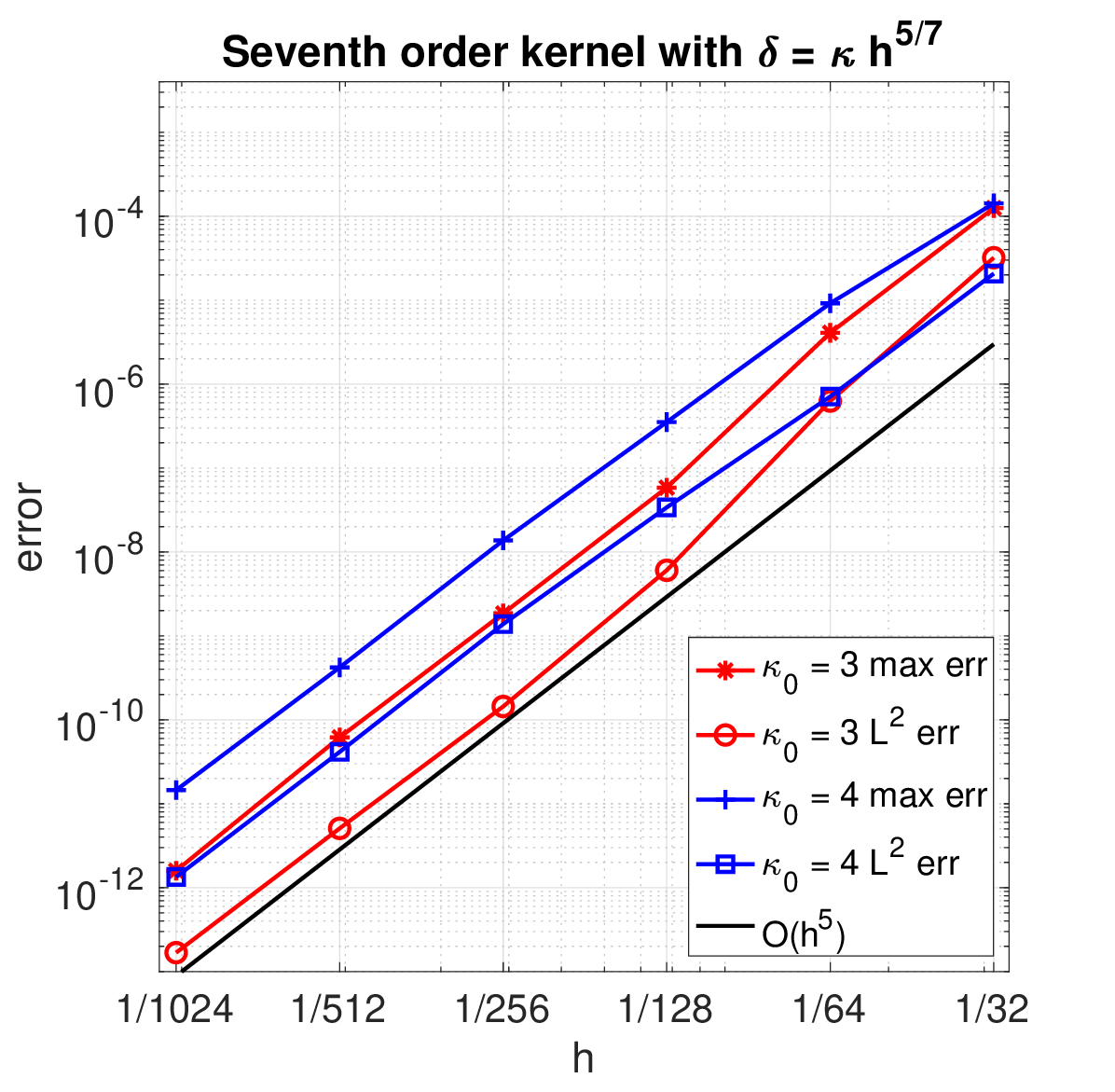}} 
\caption{Errors for the Stokes double layer on the spheroid (1,\;.5,\;.5), at grid points randomly chosen within distance $h$ to the surface. Here $\kappa$ is chosen so that $\delta/h = \kappa_0$ when $h = 1/64$.}
\label{figure:test4_stokes}
\end{figure}

The tests so far have involved the higher order versions of $s_1$ and $s_2$.  In order to test the use of $s_3$ with the
stresslet integral, we compute the integral in the identity (2.3.19) in \cite{pozbook}
\beq 
	\label{DL_identity}
	\frac{1}{8\pi} \epsilon_{jlm} \int_\Gamma x_m T_{ijk} (\bd{x_0,x}) n_k(\bd{x})dS(\bd{x}) = \chi (\bd{x}_0) \epsilon_{ilm} x_{0,m}
\eeq
where $\chi$ = 1, 1/2, 0 when $\bd{x}_0$ is inside, on, and outside the boundary.
We set $l=1$ and define $q_j(\bd{x}) = \epsilon_{j1m}x_m = (0,-x_3,x_2)$ in \eqref{stodbl}.  We use the subtracted
version \eqref{stodblsub} and regularize $T$ as in \eqref{Treg}.  Our surface is the spheroid with
semi-axes $1$, $.5$, $.5$.  We compute the solution at about $1300$ randomly selected grid points within $h$.
In Figure~\ref{figure:test4_stokes} we show results with the same $p$ and $q$ as for the molecular surface, each with two choices of $\kapz$.  The expected order is evident, especially for the larger $\kapz$.

For these examples we conclude that for orders $p = 3$, $5$, $7$ the choice
$\kapz = 2$, $3$, $4$, resp., produces reliable results and, for small enough $h$,
$\kapz = 1$, $2$, $3$ can give consistently smaller errors.


 \subsection{Extending integral values to the whole grid}  
 \label{sec:lapl_grid}

 In this section we describe and illustrate a method for obtaining
 values at all grid points from those computed as integrals at grid points near the surface.  The extension is done by a procedure due to
 A. Mayo \cite{mayo85}.  In this example we again compute single and double layer integrals for harmonic functions.
 Briefly, we compute the integrals near the surface, form a discrete Laplacian at grid points where the stencil crosses the surface,  
 extend the discrete Laplacian to be zero elsewhere, and then invert to obtain values on the entire grid.  The version we present
 leads to $O(h^4)$ accuracy for the grid values and also their first differences; the differences can be used to compute the gradient. 
 We use a fourth order discrete Laplacian described below.
 
For our test problem we use the molecular surface \eqref{molesurf}, translated to be centered at 
$\pp_0 = (1.5,1.5,1.5)$.  We embed the surface in a cube $[0, 3]^3$ and partition with a
grid of size $h$, so that there are $3/h$ points in each direction. We use a known exact solution
similar to \eqref{basicu} but nonzero in the outer region,
\beq  u_-(\yy) = (\sin{z_1} + \sin{z_2})\exp{z_3}\,,
   \quad u_+(\yy) = 1/|\zz|\,, \quad \zz = \yy - \pp_0  \eeq
Again  $ u(\yy) = \Sl(\yy) + \Dl(\yy)$ with $f = [\pa u/\pa n]$ and $g = -[u]$.

\medskip
To compute the numerical solution we proceed as follows:

\begin{enumerate}
	\item Compute $u_b(\yy) = \Sl(\yy) + \Dl(\yy)$ for grid points $\yy$ on the faces of the cube.
No special treatment is needed since the integrands are smooth.

	\item Extend the function $u_b$ on the boundary of the cube to a smooth function $w$ inside the cube.
The formula for the extension is given in Appendix A.  Set $v = u - w$, so that $v = 0$ on
the cube boundary.

	\item Compute $u_{int} = \Sl + \Dl$ at grid points within $4h$ of the surface, using seventh order
regularization with $\del = h^{5/7}$.

	\item Compute $\Delta_h u_{int}$, the discrete fourth order Laplacian of $u_{int}$,
at grid points within $2h$ of the surface.  Compute $\Delta_h w$ at all grid points.
Set $F_h = \Delta_h u_{int} - \Delta_h w $ at grid points within $2h$ and $F_h = - \Delta_h w$ elsewhere.

	\item Solve  $\Delta_h v_h = F_h$ on the cube with $v_h = 0$ on the boundary, using a discrete
sine transform version of FFT.  Set $u_h = v_h + w$, so that $u_h$ approximates $u$.

	\item If desired, compute first differences of $u_h$ such as $\left( u_h(x+h,y,z) - u_h(x,y,z) \right)/h$.
\end{enumerate}


We use a discrete Laplacian which, like the more familiar nine-point Laplacian in the plane,
is fourth order except for a correction.
This version combines the usual second order Laplacian,
using 7 points, with a sum over the corners of the cube about the central point, so that a total of 15 points are used.
The second order Laplacian is $\Delta_h^{(7)} = D_1^2 + D_2^2 +D_3^2$ where
$$ h^2\,D_1^2u(ih,jh,kh) = u((i+1)h,jh,kh) + u((i-1)h,jh,kh) - 2u(ih,jh,kh)  $$
etc.  The sum over the corner points is
$$ Cu(ih,jh,kh) = \sum_{r,s,t=-1,1} u(ih+rh,jh+sh,kh+th)  $$
and the 15-point Laplacian is
$$\Delta_h^{(15)}u = \frac{2}{3h^2}\left( h^2\Delta_h^{(7)}u \;+\; 
\frac18 Cu  \;-\; u \right)   $$
It is accurate to $O(h^4)$ except for an error $(h^2/12)\Delta^2 u$.
(See \cite{iserles}, p.169, exercise 8.11.)

\begin{figure}[!htb]
\centering
\scalebox{0.345}{\includegraphics{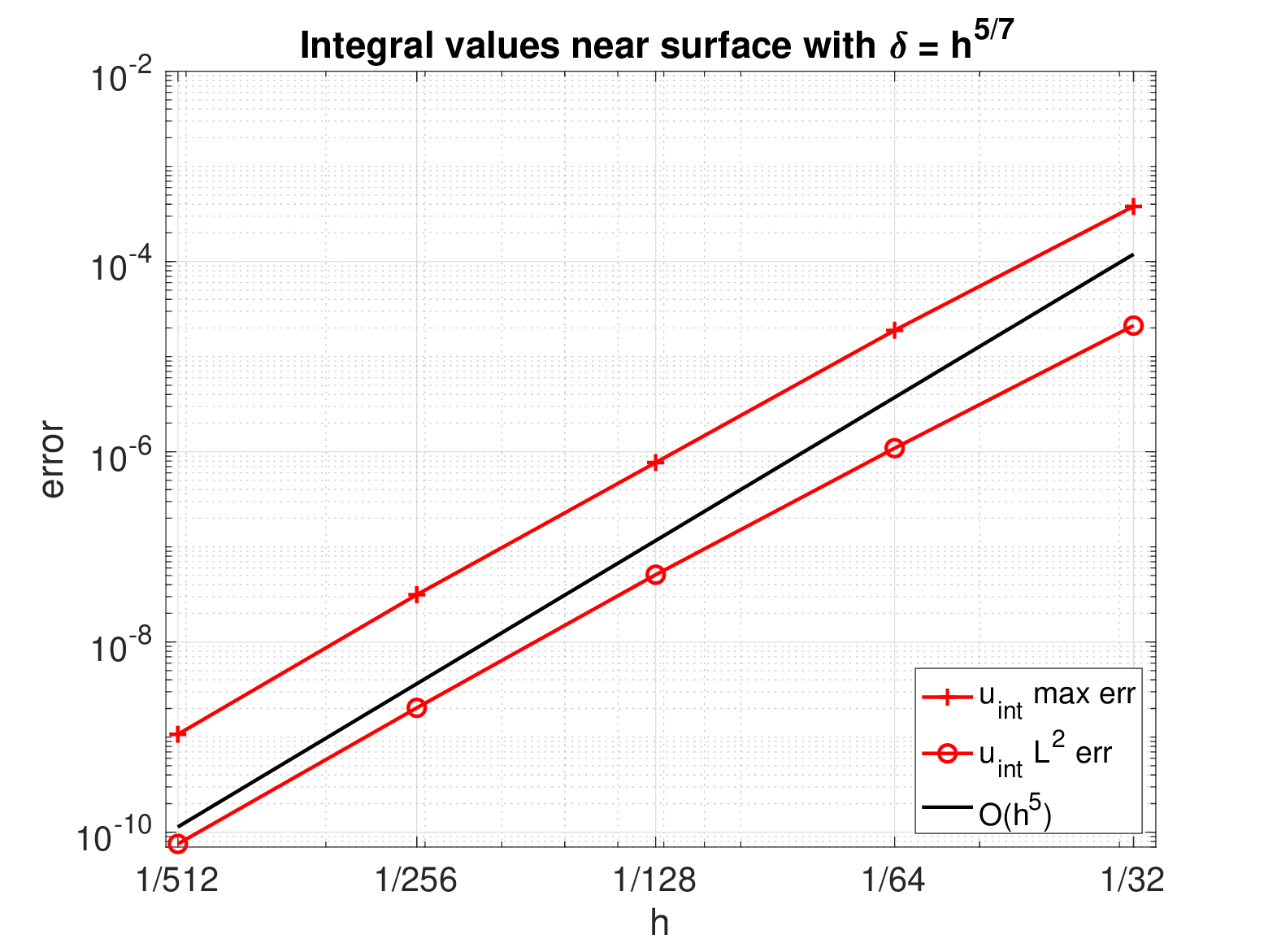}} 
\scalebox{0.345}{\includegraphics{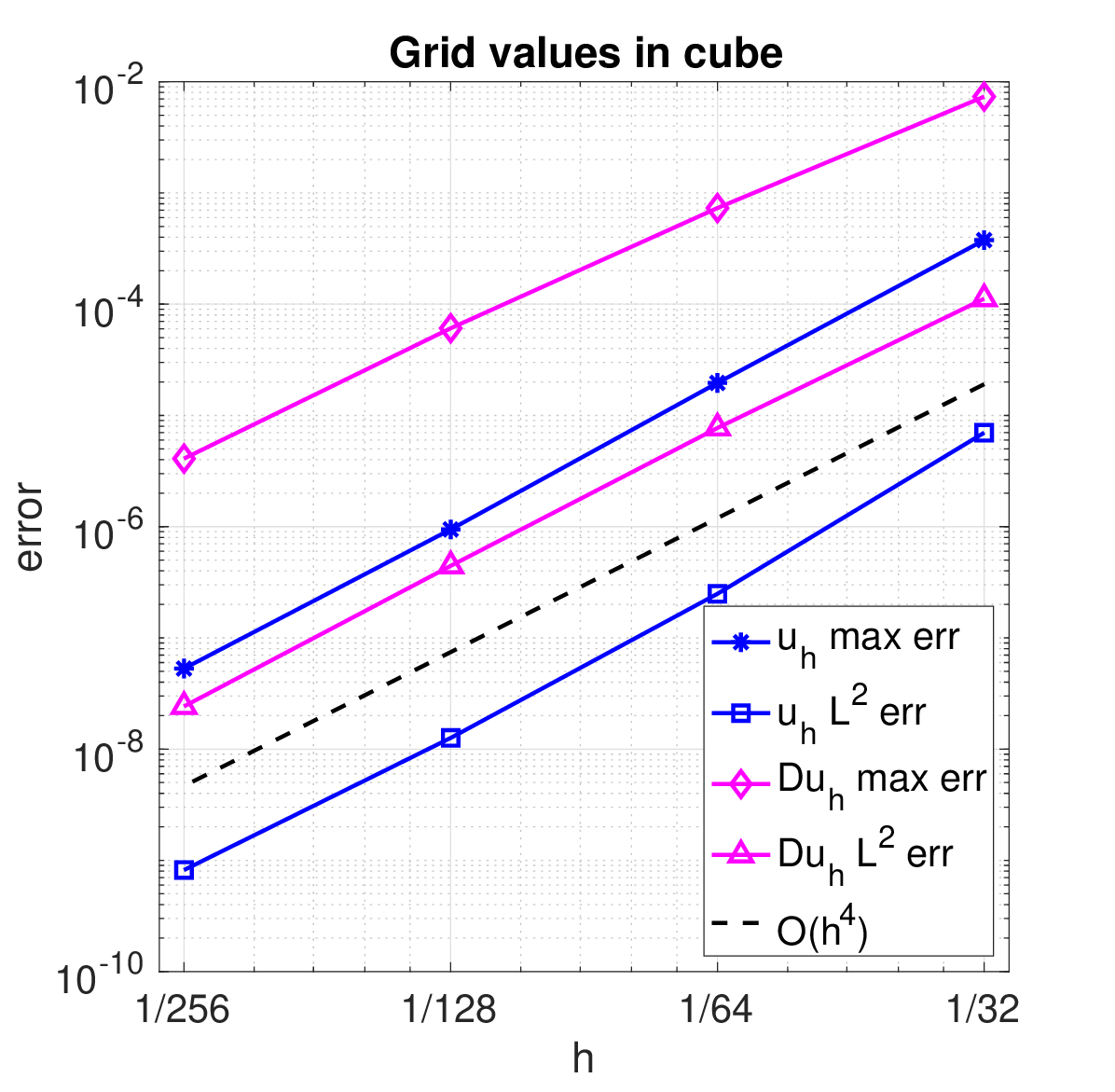}} 
\caption{Errors for the harmonic solution at grid points in $[0,3]^3$ with the molecular surface. For integrals near the surface (left), the seventh order kernel was used with $\delta=h^{5/7}$.  Errors on the whole grid in $u_h$ and its $x$-differences $Du_h$ (right), which are obtained using the fourth order Laplacian and FFT.}
\label{figure:lapl_grid}
\end{figure}

Results for the test problem are displayed in Figure~\ref{figure:lapl_grid}.  On the left, the errors are shown for the integral values $u_{int}$ found in step 3 near the surface.  On the right, errors are given for the grid values of $u_h$ obtained in step 5 above and for the horizontal first differences as in step 6.  Both show fourth order convergence.  As a consequence of the high order error in first differences, we can compute the velocity gradient in the exterior region with high accuracy:  Assuming the first differences have accuracy $O(h^4)$, we can use a fourth order accurate difference formula to calculate the gradient at grid points where the stencil does not cross the interface; the fourth order difference will be $O(h^4)$ accurate since it is a linear combination of the first differences.  We illustrate this process in Sect.~\ref{sec:stokes_grid}, where we compute the gradient of the Stokes pressure before inverting the Laplacian for the Stokes velocity.

The fourth order accuracy of the grid values and their first differences may appear surprising.
The values computed near the surface have accuracy about $O(h^5)$, so that the error in the discrete Laplacian is about $O(h^3)$.  
Because this error is restricted to an $O(h)$ neighborhood of the surface, it can be written as the discrete divergence
$\nabla_h\cdot$ of an $O(h^4)$ grid function $E$, according to Lemma 2.2 of \cite{camcos}.
The total error on the grid has the form $\nabla_h\cdot E + \tau$, where $\tau$ is the $O(h^4)$
truncation error in $\Delta_h^{(15)} u$ away from the surface.  Lemma 2.3 in \cite{camcos} provides
a regularity estimate for the usual discrete Laplacian $\Delta_h^{(7)}$.
It implies that the solution of $\Delta_h^{(7)}w = \nabla_h\cdot E + \tau$
is bounded by $\|E\|_{max} + \|\tau\|_{max}$, uniformly in $h$, and its first differences
are bounded by this quantity times $|\log{h}|$.  If the same regularity estimate
holds for $\Delta_h^{(15)}$, as we expect, it would follow that the maximum error
in grid values of $u$ computed here, and in its first differences, are $O(h^4)$ and 
$O(h^4|\log{h}|)$ respectively, consistent with the errors observed here.


\subsection{Stokes flow around a translating body}
\label{sec:stokes_grid}

We compute the  pressure and velocity in Stokes flow around a translating solid body, at grid points in the computational domain, given the force $\ff$ on the surface.  To solve the Stokes equations~\eqref{stokeseqns} we follow the procedure of the last section twice, first for the pressure and then for the velocity.  The pressure and velocity are given as integrals in~\eqref{pressure} and~\eqref{stosgl}.  We compute the integrals at grid points near the boundary $\Gamma$ of the solid as in Sect. 2.  

\begin{enumerate}
	\item Since the pressure $p$ is harmonic, the first stage is to compute $p$ just as in steps 1--5 of Sect.~\ref{sec:lapl_grid}.  To find $\Delta_h p$ in step 4, we extend $p$ to be zero inside the body.  In place of step 6, we compute the gradient $\nabla_h p$ using fourth order accurate formulas; we use a one-sided formula near the boundary of the computational box and otherwise use a standard symmetric formula. 

	\item In the second stage we compute the velocity $\uu$ by the same steps 1--5, except that $\nabla_h p$ is needed in step 4:  For grid points within $2h$ of $\Gamma$, the discrete Laplacian $F_h$ of $\vv_h = \uu_h - \ww$ is set to $\Delta_h\uu_{int} - \Delta_h \ww$ as before, whereas for other grid points we set $F_h = \nabla_h p - \Delta_h\ww$, in accordance with \eqref{stokeseqns}. As described in Sect.~\ref{sec:lapl_grid}, the gradient of $p_h$ is expected to be fourth order accurate at grid points where the stencil does not cross the interface.  In step 4 the velocity is extended to be constant inside the body.  After solving for $\uu_h$ we compute $\nabla_h\uu_h$ on the whole grid using the fourth order difference formulas.
\end{enumerate}

Our test problem is Stokes flow past a spheroid with semiaxes $(1,.5,.5)$, translating with speed
$\UU = (1,0,0)$.  Explicit formulas for the pressure, velocity, and surface force were derived in~\cite{chwang} and~\cite{liron}.  Again we compute in a cube $[0,3]^3$ with grid size $h$ and place the center of the spheroid at $\pp_0 = (1.5,1.5,1.5)$.  We use seventh order regularization with $\delta=\kappa h^{5/7}$ with $\kappa$ chosen as before such that $\delta/h=\kappa_0=4$ for $h=1/64$. We define the error in the velocity difference as $\nabla_h\uu_h - \nabla_h\uu_{exact}$, by applying the same fourth order difference to the exact velocity. We compute the error at each point in the velocity vector $\uu_h$ and in the velocity difference tensor $\nabla_h\uu_h$ using the Euclidean norm, then compute the max and $L^2$ norms as usual. 

\begin{figure}[!htb]
\centering
\scalebox{0.6}{\includegraphics{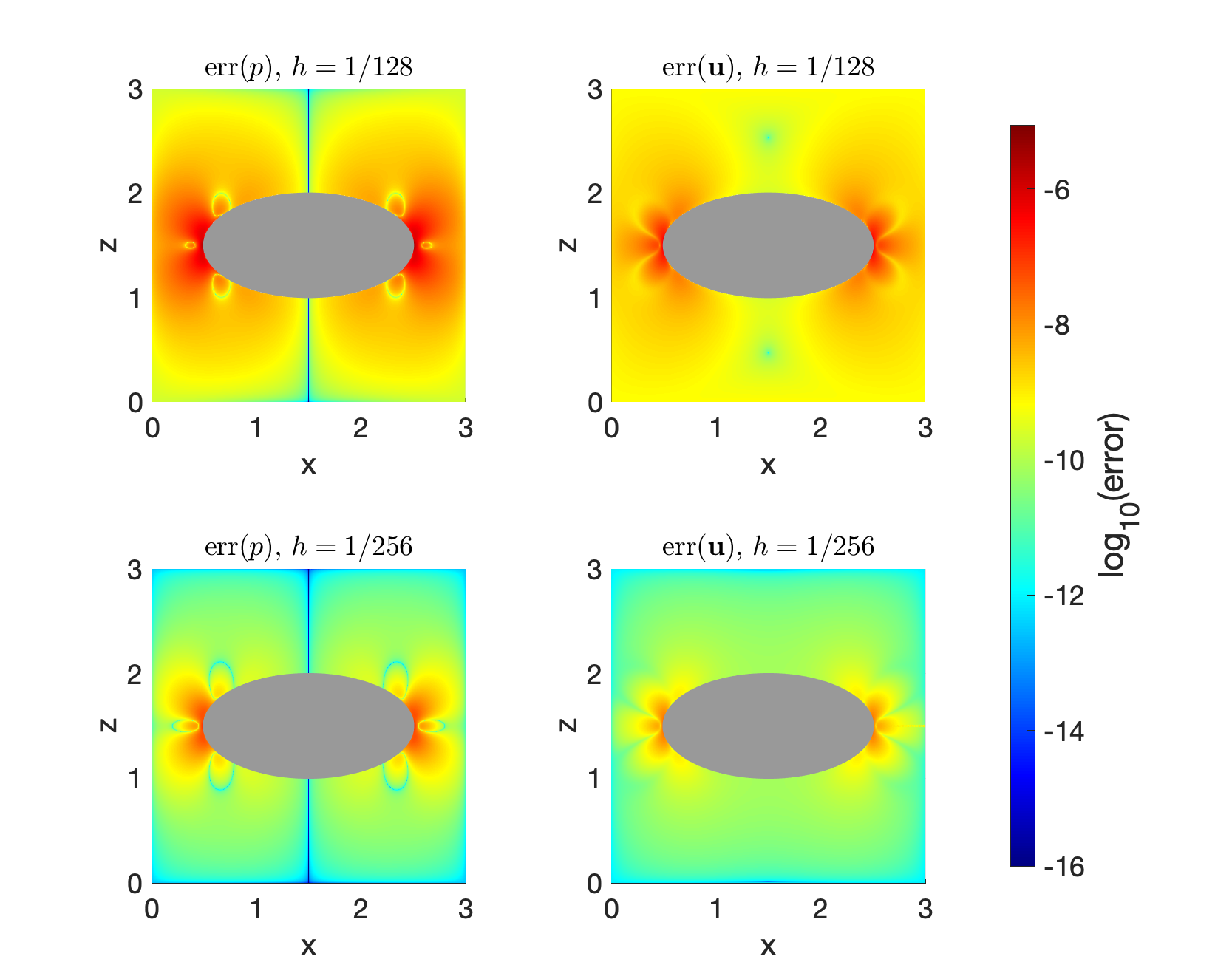}} 
\caption{Stokes solution for a spheroid $(1,.5,.5)$ translating with velocity $(1,0,0)$. Error in pressure (left) and velocity (right) in the plane $y=1.5$ for $h=1/128$ (top) and $h=1/256$ (bottom).}
\label{figure:stokes_grid_plane}
\end{figure}

\begin{figure}[!htb]
\centering
\scalebox{0.375}{\includegraphics{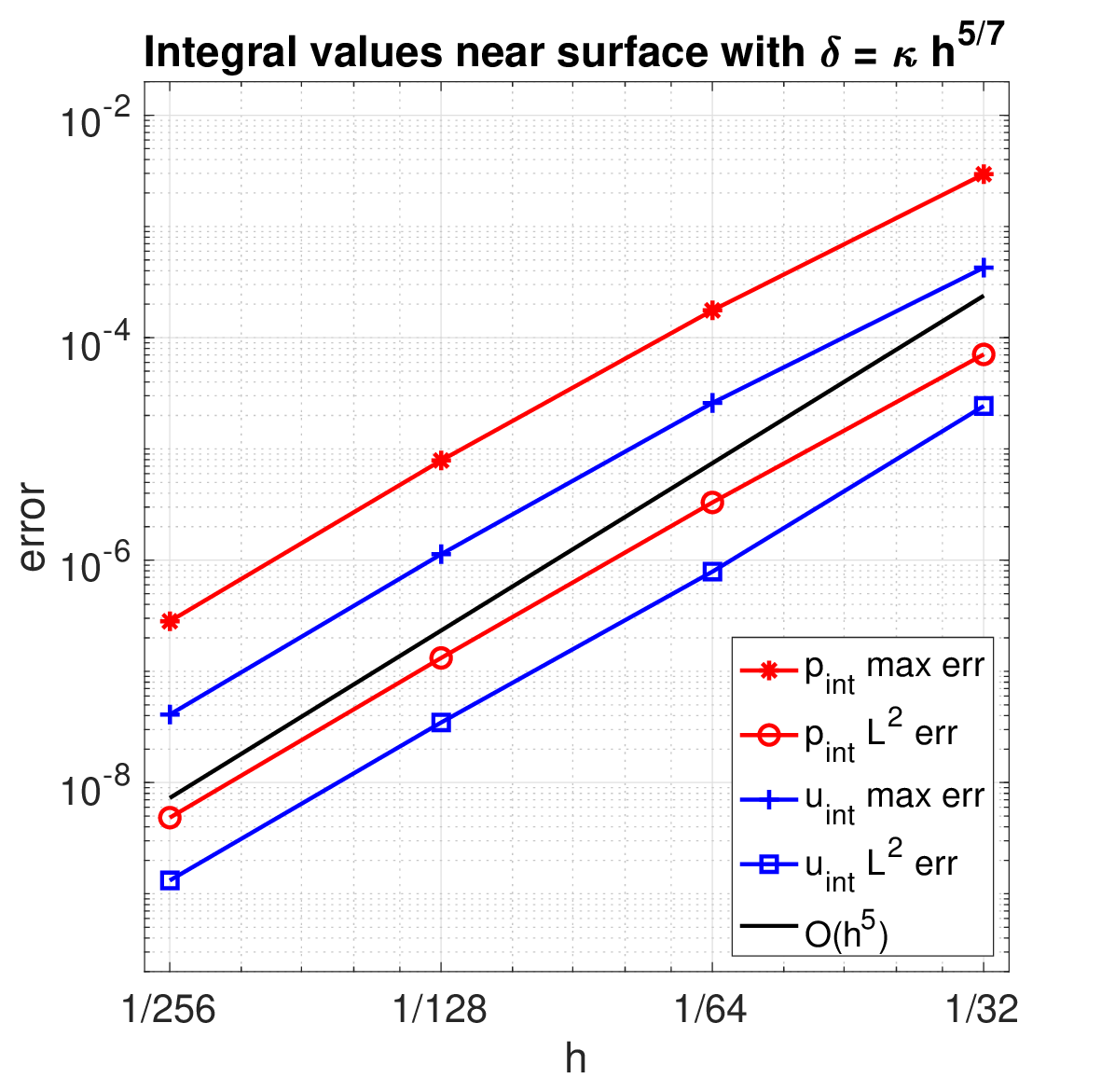}} 
\scalebox{0.375}{\includegraphics{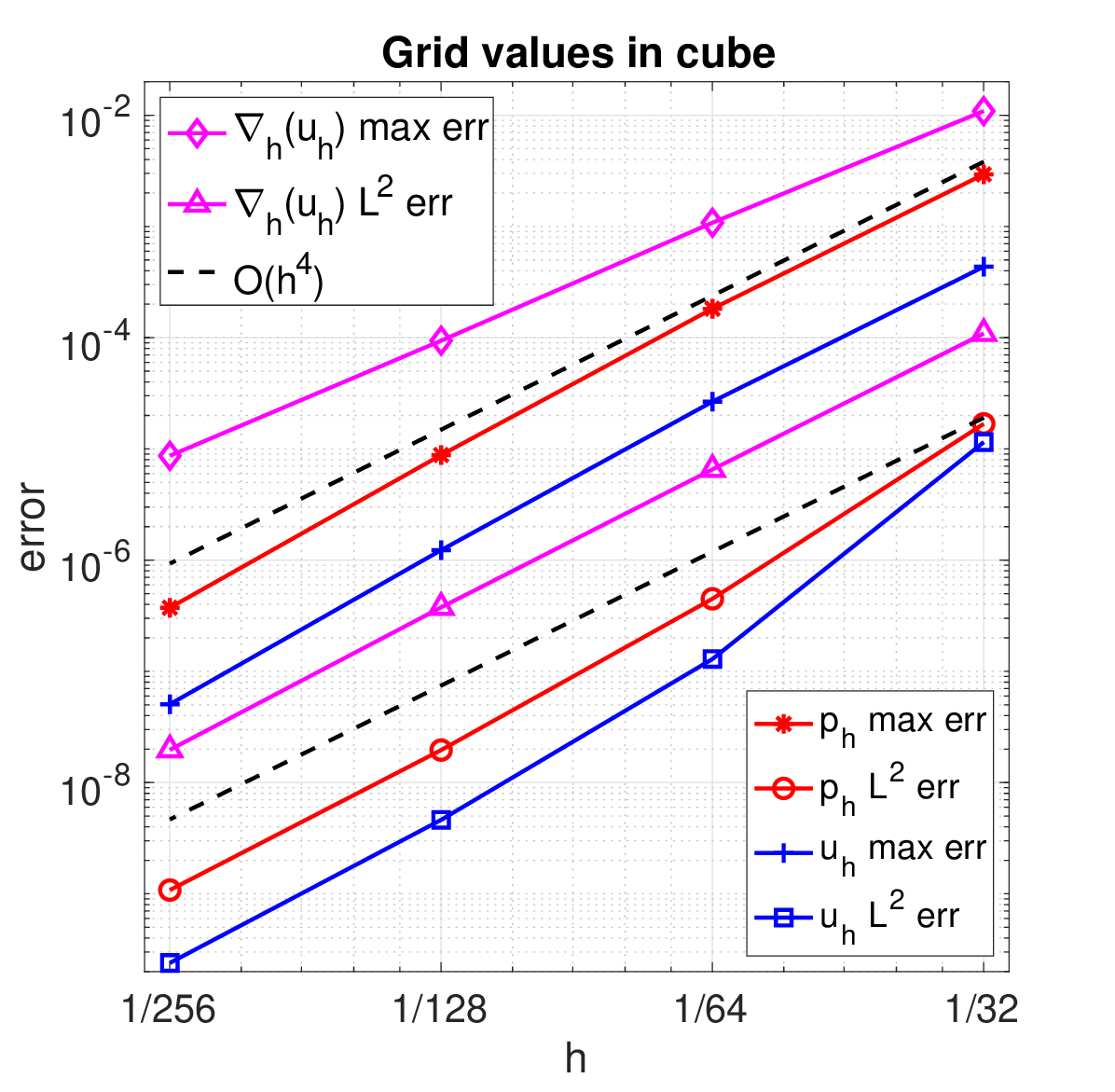}} 
\caption{Errors for the Stokes solution at grid points in $[0,3]^3$ around a translating body. For integrals near the surface (left), the seventh order kernel was used with $\delta=\kappa h^{5/7}$, $\kappa_0=4$.  Errors on the whole grid in $p_h$, $\uu_h$ and its fourth order difference tensor $\nabla_h(\uu_h)$ (right), which are obtained using the fourth order Laplacian and FFT.}
\label{figure:stokes_grid}
\end{figure}

Figure~\ref{figure:stokes_grid_plane} shows the errors in pressure (left) and velocity (right) in the plane $y=1.5$. The error is largest near the surface and decays away from it. The convergence rates are shown in Figure~\ref{figure:stokes_grid}. On the left, the errors for the integral values $p_{int}$ and $\uu_{int}$ near the surface are approaching the expected fifth order accuracy as the grid size decreases. On the right, the errors for the grid values in pressure $p_h$, velocity $\uu_h$, and differences of velocity $\nabla_h\uu_h$ are shown. The pressure and velocity are checked at grid points outside or on the boundary, and the velocity differences are checked on the whole cube. They mostly show the expected fourth order convergence or better. Since $\nabla_h\uu_{exact}$ is within $O(h^4)$ of the exact velocity gradient in the exterior region, we conclude that $\nabla_h\uu_h$ has similar accuracy.

\subsection{Surfaces close to each other in Stokes flow}
\label{sec:stokes_2surf}

Nearly singular integrals can be needed when surfaces are close to each other. We consider two interfaces $\Gamma_b, b=1,2$ between fluids with different viscosities that are near touching. The interface velocity satisfies an integral equation~\cite{pozbook}
\begin{align}
	\label{IE-2Surf}
	(\lambda_b+1) u_i(\textbf{x}_0) = &-\frac{1}{4\pi \mu_0}\sum_{m=1}^{2}\int_{\Gamma_m} S_{ij}(\textbf{x}_0,\textbf{x}) [f]_j(\textbf{x})dS(\textbf{x}) \nonumber\\
	&+ \sum_{m=1}^{2}\frac{\lambda_m-1}{4\pi}\int_{\Gamma_m} T_{ijk} (\textbf{x}_0,\textbf{x}) u_j(\textbf{x}) n_k(\textbf{x})dS(\textbf{x}), \quad \textbf{x}_0\in \Gamma_b
\end{align}
where $u_i$ is the $i$-th velocity component. We take the external viscosity to be $\mu_0=1$ and internal viscosities $\mu_1 = \mu_2 = 2$, with the viscosity ratio for each interface $\lambda_b=\mu_b/\mu_0=2$. While the velocity across the interface is continuous, the surface force has a discontinuity given by $[\bd{f}] = 2\gamma H \bd{n} - \nabla_S\gamma$, where $\gamma$ is the surface tension, $H$ is the mean curvature, $\bd{n}$ is the outward unit normal, and $\nabla_S = \nabla - \nn(\nn\cdot\nabla)$ is the surface gradient.

For our test we choose two spheroids with semiaxes $(1,.5,.5)$, one centered at $(0,0,-.5-\epsilon)$, and the other centered at $(0,0,.5)$ and rotated by $30^\circ$ in the $xy$-plane. We take the distance between the surfaces to be $\epsilon = 1/16^3$, and define the surface tension on each spheroid as $\gamma = 1+(z-z_c)^2$, where $z_c$ is the $z$-coordinate of the center of the spheroid.  

To solve the integral equation~\eqref{IE-2Surf}, we use GMRES with a tolerance of $10^{-10}$, which takes about 12 iterations in each of our runs. Initially, two Stokeslet integrals are evaluated for each surface. For each target point on a given surface, one integral is evaluated on the same surface and we use formulas of Sect.~\ref{sec:onSurf}. The second integral is evaluated over the other surface, and we use formulas of Sect.~\ref{sec:offSurf} when the target is near the other surface, with $b/\delta \leq 8$. Otherwise, we evaluate the Stokeslet integral~\eqref{stosgl} without special treatment. At each iteration, a similar procedure is followed for the stresslet integrals. We use a treecode~\cite{wang-krasny-tlupova-20} for efficient summation in the far field, using parameters as suggested in~\cite{siebor-tlupova}. The exact solution is not known in this case, and the error at grid size $h$ is checked empirically using the solutions at $h$ and $h/2$. 

\begin{figure}[!htb]
\centering
\scalebox{0.375}{\includegraphics{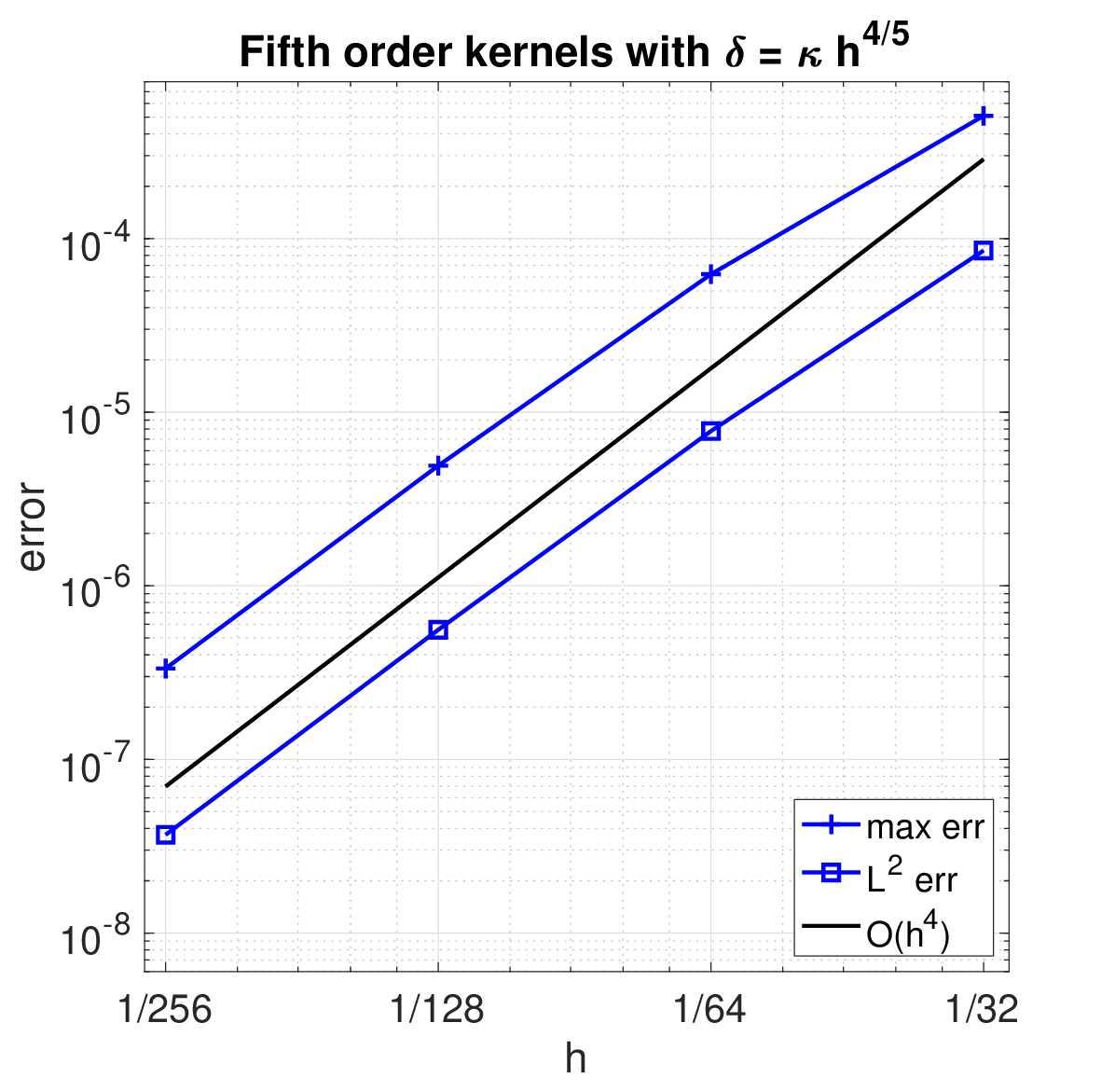}} 
\scalebox{0.375}{\includegraphics{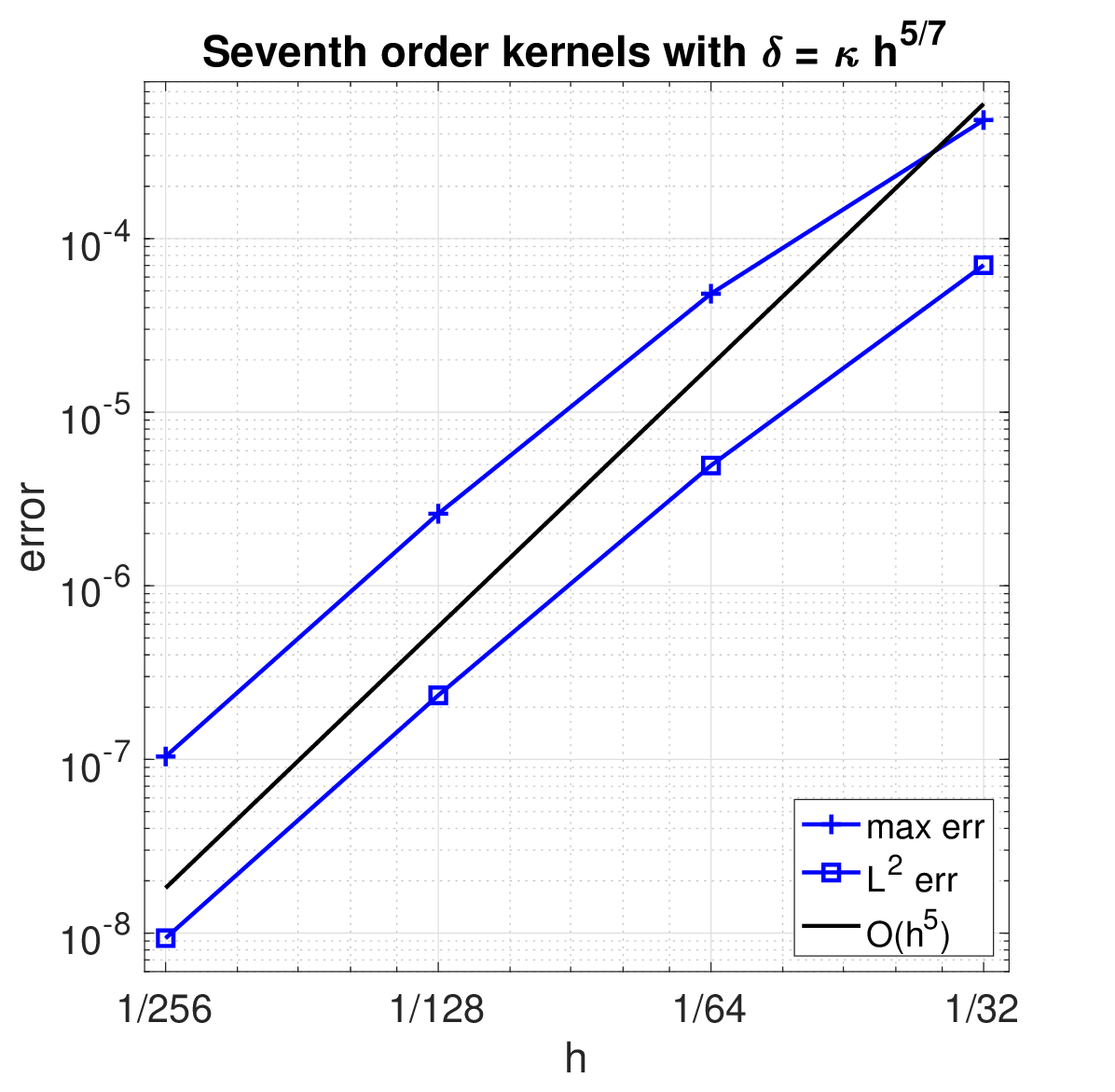}} 
\caption{Errors in the interface velocity on the two spheroids separated by $\epsilon=1/16^3$. Here $\kappa$ is chosen so that $\delta/h = \kappa_0$ when $h=1/64$. We use $\kappa_0=3$ with the fifth order kernels and $\kappa_0=4$ with the seventh order kernels.}
\label{figure:2surf_errors}
\end{figure}

Figure~\ref{figure:2surf_errors} shows the errors in the velocity of both interfaces using regularized kernels of fifth order with $\delta=\kappa h^{4/5}$ and seventh order with $\delta=\kappa h^{5/7}$. Using the guidance of Sect. 4.2, we again choose $\kappa_0=3$ and $\kappa_0=4$ for the two orders, respectively. The errors are near identical when considered for individual spheroids. For the fifth order regularization, the order of convergence is evident. For the seventh order regularization, the expected order of convergence becomes evident as the grid size gets reduced. This is similar to the results in Figure~\ref{figure:lapl_grid}, left, where the single and double layer potentials are combined for the harmonic solution, as well as Figure~\ref{figure:stokes_grid}, left, for the Stokes pressure and velocity. Results presented here are for spheroids that are within a small distance $\epsilon=1/16^3$ from each other. When using a larger distance $\epsilon=1/16$, the overall errors are almost identical. 

\begin{figure}[!htb]
\centering
\scalebox{0.6}{\includegraphics{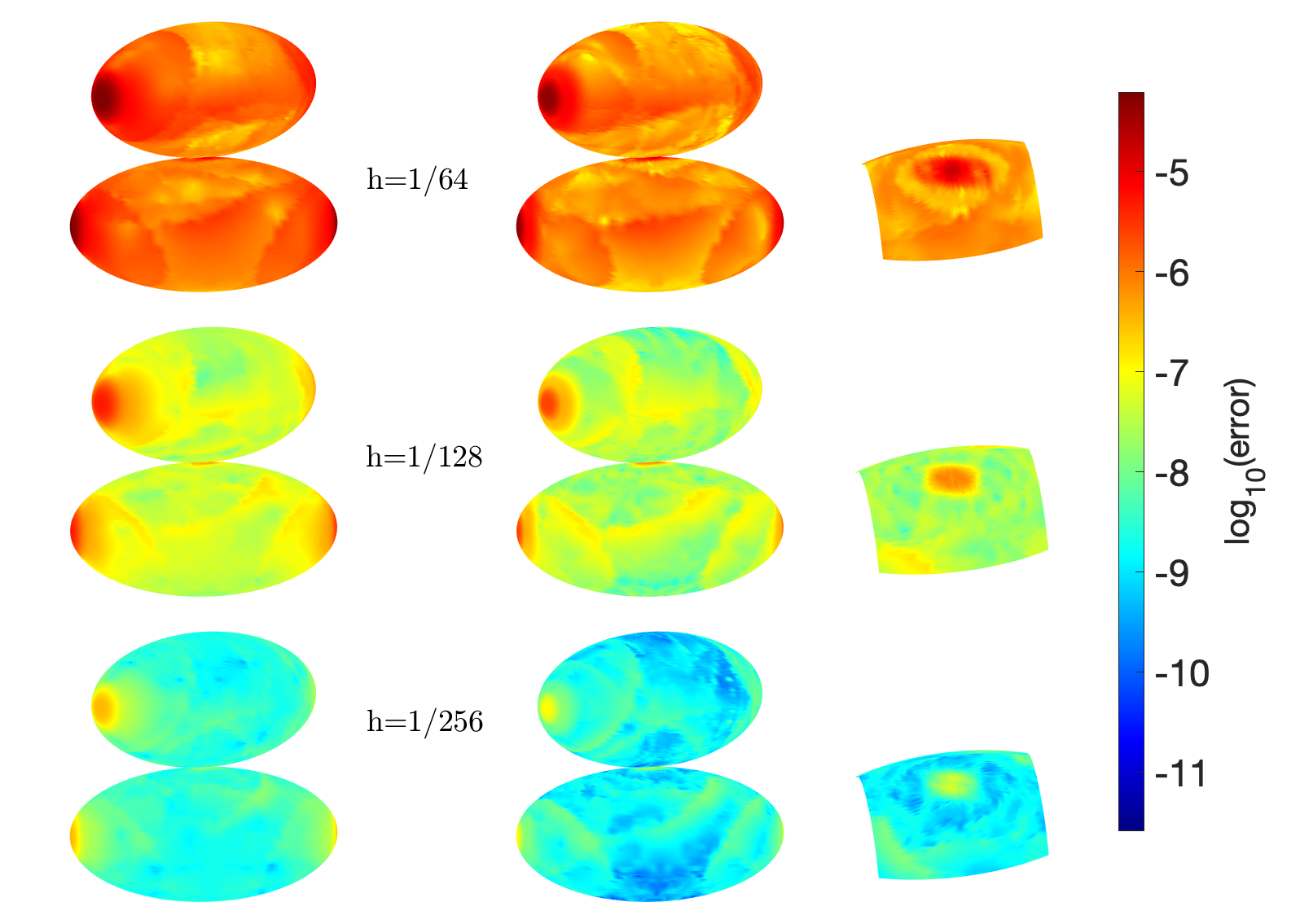}} 
\caption{Velocity errors on the interfaces on the two spheroids separated by $\epsilon=1/16^3$. Left: fifth order kernels with $\delta=\kappa h^{4/5}$, $\kappa_0=3$. Middle: seventh order kernels with $\delta=\kappa h^{5/7}$, $\kappa_0=4$. Right: near touching points on the lower surface using the seventh order kernels.}
\label{figure:2surf}
\end{figure}

The surface distributions of the errors for three different grid sizes are shown in Figure~\ref{figure:2surf}. The errors using the fifth order kernels are shown in the left column, and the errors using the seventh order kernels are shown in the middle column. The right column shows the errors at points on the lower surface that are near touching the upper surface. The maximum error over the near touching points and the maximum error over the points on the edges of the spheroid are within a factor of two of each other.

\section{Conclusions}
We have derived formulas for regularizing single or double layer integrals on surfaces, for harmonic functions or Stokes flow.  The error is high order in the regularization parameter $\delta$.  These formulas can be used directly to compute the integrals for target points on or near the surface.  The significance of the high order regularization is that it allows $\delta$ to be large enough to control the discretization error with a conventional quadrature.  The derivation is based on the analysis in~\cite{extrap}.  We illustrate the use of these formulas to compute values at grid points near the surface and then extend to values on the entire grid, as well as computing integrals on surfaces that are close together.

In our examples we choose $\delta = \kappa h^q$, $q \leq 1$, where $h$ is the mesh spacing.  There cannot be a definite rule for choosing the relationship between $\delta$ and $h$, but we have developed a guideline that has produced the expected orders of accuracy in our tests.  The regularization error is $O(\delta^p)$ with $p = 3$, $5$, or $7$, so that the predicted error is $O(h^{pq})$.  In our examples we take $q = 2/3$, $4/5$, and $5/7$ for $p = 3, 5, 7$.
To choose $\kappa$, we first choose a relatively large $h_0$ which resolves the surface, in our case $h_0 = 1/64$.  We set $\delta_0 = \kappa_0 h_0$.
We have taken $\kappa_0 = 2$, $3$, $4$ for $p = 3, 5, 7$, although $\kappa_0 = 1$, $2$, $3$ may give smaller errors.  Finally we choose $\kappa$ so that $\delta_0 = \kappa h_0^q$.

In \cite{extrap} we computed the integrals using a simple regularization with several choices of $\delta$ and extrapolated to obtain high order.  The present method is similar in that it is based on the same analysis.  It is simpler in that it uses only one $\delta$, but the formulas are more involved.  It gives the predicted order more reliably when $p = 7$. 

For a surface layer potential defining a harmonic function $u$, once values of $u$ are computed at grid points near the surface, the solution can be found on the entire grid using a procedure of~\cite{mayo85}.  We form the discrete Laplacian $\Delta_h u$ at grid points where the stencil crosses the surface.  We set $\Delta_h u = 0$ at other grid points and then invert $\Delta_h$ using fft's to obtain grid values of $u$.  We compute the integrals with $p = 7$ and $q= 5/7$ and use a fourth order version of $\Delta_h$, obtaining grid values that are accurate to $O(h^4)$; see Sect. 4.3.  We use a similar procedure for Stokes flow in Sect. 4.4, solving first for the pressure and then for the velocity.  


\begin{appendix}
\section{Extending a function from the surface of a cube to the interior}
We begin with the simpler task of extending a function into a square $[0,1]^2$ from the four edges. 
Assume $g$ is a smooth function given on the boundary; in particular, we assume that the values 
at a corner along the two edges are equal.
The extension is defined as $\gbar(x_1,x_2) = E(x_1,x_2) - C(x_1,x_2)$
where
\beq E(x_1,x_2) = x_1g(1,x_2) + (1-x_1)g(0,x_2) + x_2g(x_1,1) + (1-x_2)g(x_1,0)  \eeq
\beq C(x_1,x_2) = x_1x_2g(1,1) + x_1(1-x_2)g(1,0) + (1-x_1)x_2g(0,1) + (1-x_1)(1-x_2)g(0,0) \eeq
The function $E$ performs linear interpolation in each direction, and $C$ makes corrections using corner values.
This formula was provided by Thomas Witelski, and it motivates the following formula for a cube.

Now suppose $g(x_1,x_2,x_3)$ is a smooth function on the faces of the cube $[0,1]^3$.
We define an extension in the form $\gbar = F - E + C$, starting with
linear interpolation from the faces,
\begin{multline}
 F(x_1,x_2,x_3) = (1-x_1)g(0,x_2,x_3) + x_1g(1,x_2,x_3) + (1-x_2)g(x_1,0,x_3) \\
    + x_2g(x_1,1,x_3) + (1-x_3)g(x_1,x_2,0) + x_3g(x_1,x_2,1) 
\end{multline}
We will write this and two further sums in a more compact form.  Define
\beq \ell(w,p) = 1-w\; \mbox{if} \;p=0\;\; \mbox{and} \;\;\ell(w,p) = w\; \mbox{if} \;p=1 \eeq
Then
\beq F(x_1,x_2,x_3) = \sum_{i=1}^3 \sum_{r=0}^1 \ell(x_i,r)g(z_1,z_2,z_3) \eeq
where $z_i = r$ and $z_j = x_j$, $z_k = x_k$ for $j,k \neq i$.
The function $E$ uses values on the edges,
\beq E(x_1,x_2,x_3) = \sum_{i=1}^3 \sum_{r=0}^1 \sum_{s=0}^1 \ell(x_j,r)\ell(x_k,s) g(z_1,z_2,z_3) \eeq
Here for each $i$, if $j,k$ are the other two indices, then $z_i = x_i$, $z_j = r$, and $z_k = s$.
For example if $i=1$, $r=0$, $s = 1$, then  $j=2$, $k=3$,
and the term is $(1-x_2)x_3 g(x_1,0,1)$.
Finally we need another sum that matches the values of $g$ at the corners,
\beq C(x_1,x_2,x_3) = \sum_{r=0}^1 \sum_{s=0}^1 \sum_{t=0}^1  
     \ell(x_1,r)\ell(x_2,s)\ell(x_3,t) g(r,s,t)  \eeq
The function $\gbar = F - E + C$ defined inside the cube agrees with $g$ on the faces.
 \end{appendix}

\section*{Declarations}

\subsection*{Conflict of interest}

The authors declare no competing interests.

\section*{Acknowledgment}
We thank Thomas Witelski for providing the extension formula for a square.


\bibliographystyle{plain}
\bibliography{hireg}

\end{document}